\documentclass[11pt]{article}
\usepackage{amsmath, upgreek}
\usepackage{amssymb}
\usepackage{amsfonts}
\usepackage{amsmath,tabu}
\usepackage{cite}
\usepackage{color}

\usepackage{xcolor}
\usepackage{amscd}
\usepackage{xspace}
\usepackage{verbatim}
\usepackage{graphicx}
\newcommand{\mysection}[1]{
\section{#1}\setcounter{equation}{0}}
\title{\bf Isolated singularities of solutions of a 2-D  diffusion equation with mixed reaction
}
\author{{\bf Yimei Li\footnote{\noindent (1)School of Mathematics and Statistics, Beijing Jiaotong University, Beijing 100044, PR China.  (2) Universit\'e Sorbonne Paris Nord, CNRS UMR 7539, LAGA, 93430 Villetaneuse, France.. E-mail: liyimei@bjtu.edu.cn}} \\
 {\bf Laurent V\'eron \footnote{\noindent
Institut Denis Poisson, CNRS UMR 7013,  Universit\'e de Tours, 37200 Tours, France. 
          E-mail: veronl@univ-tours.fr}}\\[2mm]
}
 
\date{}
\begin{document}
 \maketitle


\newcommand{\txt}[1]{\;\text{ #1 }\;}
\newcommand{\tbf}{\textbf}
\newcommand{\tit}{\textit}
\newcommand{\tsc}{\textsc}
\newcommand{\trm}{\textrm}
\newcommand{\mbf}{\mathbf}
\newcommand{\mrm}{\mathrm}
\newcommand{\bsym}{\boldsymbol}
\newcommand{\scs}{\scriptstyle}
\newcommand{\sss}{\scriptscriptstyle}
\newcommand{\txts}{\textstyle}
\newcommand{\dsps}{\displaystyle}
\newcommand{\fnz}{\footnotesize}
\newcommand{\scz}{\scriptsize}
\newcommand{\be}{\begin{equation}}
\newcommand{\bel}[1]{\begin{equation}\label{#1}}
\newcommand{\ee}{\end{equation}}
\newcommand{\eqnl}[2]{\begin{equation}\label{#1}{#2}\end{equation}}
\newcommand{\barr}{\begin{eqnarray}}
\newcommand{\earr}{\end{eqnarray}}
\newcommand{\bars}{\begin{eqnarray*}}
\newcommand{\ears}{\end{eqnarray*}}
\newcommand{\nnu}{\nonumber \\}
\newtheorem{subn}{\name}
\renewcommand{\thesubn}{}
\newcommand{\bsn}[1]{\def\name{#1}\begin{subn}}
\newcommand{\esn}{\end{subn}}
\newtheorem{sub}{\name}[section]
\newcommand{\dn}[1]{\def\name{#1}}   
\newcommand{\bs}{\begin{sub}}
\newcommand{\es}{\end{sub}}
\newcommand{\bsl}[1]{\begin{sub}\label{#1}}
\newcommand{\bth}[1]{\def\name{Theorem}
\begin{sub}\label{t:#1}}
\newcommand{\blemma}[1]{\def\name{Lemma}
\begin{sub}\label{l:#1}}
\newcommand{\bcor}[1]{\def\name{Corollary}
\begin{sub}\label{c:#1}}
\newcommand{\bdef}[1]{\def\name{Definition}
\begin{sub}\label{d:#1}}
\newcommand{\bprop}[1]{\def\name{Proposition}
\begin{sub}\label{p:#1}}

\newcommand{\R}{\eqref}
\newcommand{\rth}[1]{Theorem~\ref{t:#1}}
\newcommand{\rlemma}[1]{Lemma~\ref{l:#1}}
\newcommand{\rcor}[1]{Corollary~\ref{c:#1}}
\newcommand{\rdef}[1]{Definition~\ref{d:#1}}
\newcommand{\rprop}[1]{Proposition~\ref{p:#1}}
\newcommand{\BA}{\begin{array}}
\newcommand{\EA}{\end{array}}
\newcommand{\BAN}{\renewcommand{\arraystretch}{1.2}
\setlength{\arraycolsep}{2pt}\begin{array}}
\newcommand{\BAV}[2]{\renewcommand{\arraystretch}{#1}
\setlength{\arraycolsep}{#2}\begin{array}}
\newcommand{\BSA}{\begin{subarray}}
\newcommand{\ESA}{\end{subarray}}
\newcommand{\BAL}{\begin{aligned}}
\newcommand{\EAL}{\end{aligned}}
\newcommand{\BALG}{\begin{alignat}}
\newcommand{\EALG}{\end{alignat}}
\newcommand{\BALGN}{\begin{alignat*}}
\newcommand{\EALGN}{\end{alignat*}}
\newcommand{\note}[1]{\textit{#1.}\hspace{2mm}}
\newcommand{\Proof}{\note{Proof}}
\newcommand{\qeda}{\hspace{10mm}\hfill $\square$}
\newcommand{\qed}{\\
${}$ \hfill $\square$}
\newcommand{\Remark}{\note{Remark}}
\newcommand{\modin}{$\,$\\[-4mm] \indent}
\newcommand{\forevery}{\quad \forall}
\newcommand{\set}[1]{\{#1\}}
\newcommand{\setdef}[2]{\{\,#1:\,#2\,\}}
\newcommand{\setm}[2]{\{\,#1\mid #2\,\}}
\newcommand{\mt}{\mapsto}
\newcommand{\lra}{\longrightarrow}
\newcommand{\lla}{\longleftarrow}
\newcommand{\llra}{\longleftrightarrow}
\newcommand{\Lra}{\Longrightarrow}
\newcommand{\Lla}{\Longleftarrow}
\newcommand{\Llra}{\Longleftrightarrow}
\newcommand{\warrow}{\rightharpoonup}
\newcommand{
\paran}[1]{\left (#1 \right )}
\newcommand{\sqbr}[1]{\left [#1 \right ]}
\newcommand{\curlybr}[1]{\left \{#1 \right \}}
\newcommand{\abs}[1]{\left |#1\right |}
\newcommand{\norm}[1]{\left \|#1\right \|}
\newcommand{
\paranb}[1]{\big (#1 \big )}
\newcommand{\lsqbrb}[1]{\big [#1 \big ]}
\newcommand{\lcurlybrb}[1]{\big \{#1 \big \}}
\newcommand{\absb}[1]{\big |#1\big |}
\newcommand{\normb}[1]{\big \|#1\big \|}
\newcommand{
\paranB}[1]{\Big (#1 \Big )}
\newcommand{\absB}[1]{\Big |#1\Big |}
\newcommand{\normB}[1]{\Big \|#1\Big \|}
\newcommand{\produal}[1]{\langle #1 \rangle}

\newcommand{\thkl}{\rule[-.5mm]{.3mm}{3mm}}
\newcommand{\thknorm}[1]{\thkl #1 \thkl\,}
\newcommand{\trinorm}[1]{|\!|\!| #1 |\!|\!|\,}
\newcommand{\bang}[1]{\langle #1 \rangle}
\def\angb<#1>{\langle #1 \rangle}
\newcommand{\vstrut}[1]{\rule{0mm}{#1}}
\newcommand{\rec}[1]{\frac{1}{#1}}
\newcommand{\opname}[1]{\mbox{\rm #1}\,}
\newcommand{\supp}{\opname{supp}}
\newcommand{\dist}{\opname{dist}}
\newcommand{\myfrac}[2]{{\displaystyle \frac{#1}{#2} }}
\newcommand{\myint}[2]{{\displaystyle \int_{#1}^{#2}}}
\newcommand{\mysum}[2]{{\displaystyle \sum_{#1}^{#2}}}
\newcommand {\dint}{{\displaystyle \myint\!\!\myint}}
\newcommand{\q}{\quad}
\newcommand{\qq}{\qquad}
\newcommand{\hsp}[1]{\hspace{#1mm}}
\newcommand{\vsp}[1]{\vspace{#1mm}}
\newcommand{\ity}{\infty}
\newcommand{\prt}{\partial}
\newcommand{\sms}{\setminus}
\newcommand{\ems}{\emptyset}
\newcommand{\ti}{\times}
\newcommand{\pr}{^\prime}
\newcommand{\ppr}{^{\prime\prime}}
\newcommand{\tl}{\tilde}
\newcommand{\sbs}{\subset}
\newcommand{\sbeq}{\subseteq}
\newcommand{\nind}{\noindent}
\newcommand{\ind}{\indent}
\newcommand{\ovl}{\overline}
\newcommand{\unl}{\underline}
\newcommand{\nin}{\not\in}
\newcommand{\pfrac}[2]{\genfrac{(}{)}{}{}{#1}{#2}}

\def\ga{\alpha}     \def\gb{\beta}       \def\gg{\gamma}
\def\gc{\chi}       \def\gd{\delta}      \def\ge{\epsilon}
\def\gth{\theta}                         \def\vge{\varepsilon}
\def\gf{\phi}       \def\vgf{\varphi}    \def\gh{\eta}
\def\gi{\iota}      \def\gk{\kappa}      \def\gl{\lambda}
\def\gm{\mu}        \def\gn{\nu}         \def\gp{\pi}
\def\vgp{\varpi}    \def\gr{\rho}        \def\vgr{\varrho}
\def\gs{\sigma}     \def\vgs{\varsigma}  \def\gt{\tau}
\def\gu{\upsilon}   \def\gv{\vartheta}   \def\gw{\omega}
\def\gx{\xi}        \def\gy{\psi}        \def\gz{\zeta}
\def\Gg{\Gamma}     \def\Gd{\Delta}      \def\Gf{\Phi}
\def\Gth{\Theta}
\def\Gl{\Lambda}    \def\Gs{\Sigma}      \def\Gp{\Pi}
\def\Gw{\Omega}     \def\Gx{\Xi}         \def\Gy{\Psi}

\def\CS{{\mathcal S}}   \def\CM{{\mathcal M}}   \def\CN{{\mathcal N}}
\def\CR{{\mathcal R}}   \def\CO{{\mathcal O}}   \def\CP{{\mathcal P}}
\def\CA{{\mathcal A}}   \def\CB{{\mathcal B}}   \def\CC{{\mathcal C}}
\def\CD{{\mathcal D}}   \def\CE{{\mathcal E}}   \def\CF{{\mathcal F}}
\def\CG{{\mathcal G}}   \def\CH{{\mathcal H}}   \def\CI{{\mathcal I}}
\def\CJ{{\mathcal J}}   \def\CK{{\mathcal K}}   \def\CL{{\mathcal L}}
\def\CT{{\mathcal T}}   \def\CU{{\mathcal U}}   \def\CV{{\mathcal V}}
\def\CZ{{\mathcal Z}}   \def\CX{{\mathcal X}}   \def\CY{{\mathcal Y}}
\def\CW{{\mathcal W}} \def\CQ{{\mathcal Q}}
\def\BBA {\mathbb A}   \def\BBb {\mathbb B}    \def\BBC {\mathbb C}
\def\BBD {\mathbb D}   \def\BBE {\mathbb E}    \def\BBF {\mathbb F}
\def\BBG {\mathbb G}   \def\BBH {\mathbb H}    \def\BBI {\mathbb I}
\def\BBJ {\mathbb J}   \def\BBK {\mathbb K}    \def\BBL {\mathbb L}
\def\BBM {\mathbb M}   \def\BBN {\mathbb N}    \def\BBO {\mathbb O}
\def\BBP {\mathbb P}   \def\BBR {\mathbb R}    \def\BBS {\mathbb S}
\def\BBT {\mathbb T}   \def\BBU {\mathbb U}    \def\BBV {\mathbb V}
\def\BBW {\mathbb W}   \def\BBX {\mathbb X}    \def\BBY {\mathbb Y}
\def\BBZ {\mathbb Z}

\def\GTA {\mathfrak A}   \def\GTB {\mathfrak B}    \def\GTC {\mathfrak C}
\def\GTD {\mathfrak D}   \def\GTE {\mathfrak E}    \def\GTF {\mathfrak F}
\def\GTG {\mathfrak G}   \def\GTH {\mathfrak H}    \def\GTI {\mathfrak I}
\def\GTJ {\mathfrak J}   \def\GTK {\mathfrak K}    \def\GTL {\mathfrak L}
\def\GTM {\mathfrak M}   \def\GTN {\mathfrak N}    \def\GTO {\mathfrak O}
\def\GTP {\mathfrak P}   \def\GTR {\mathfrak R}    \def\GTS {\mathfrak S}
\def\GTT {\mathfrak T}   \def\GTU {\mathfrak U}    \def\GTV {\mathfrak V}
\def\GTW {\mathfrak W}   \def\GTX {\mathfrak X}    \def\GTY {\mathfrak Y}
\def\GTZ {\mathfrak Z}   \def\GTQ {\mathfrak Q}

\font\Sym= msam10 
\def\SYM#1{\hbox{\Sym #1}}
\newcommand{\bdw}{\prt\Gw\xspace}
\date{}
\maketitle\medskip

\noindent{\small {\bf Abstract} We study the local properties of positive solutions of the equation $-\Gd u+ae^{bu}=m\abs{\nabla u}^q$ in a punctured domain $\Gw\setminus\{0\}$ of $\BBR^2$ where $m,a,b$ are positive parameters  and $q>1$. We study particularly the existence of solutions with an isolated singularity and the local behaviour of such singular solutions.
}\medskip

\tableofcontents
\vspace{1mm}
\hspace{.05in}
\medskip

\mysection{Introduction}

The following type of nonlinear equation  
\bel{I-0}
-\Gd u +f(u) -g(|\nabla u|)=0
\ee
in a $N$ dimensional domain when $f$ and $g$ are nondecreasing continuous functions is an interesting model for studying the interaction between the diffusion, the reaction and the absorption
 when the last two quantities are of a totally different nature. The most interesting case correspond to the case where the two reaction functions $f$ and $g$  increasing in their respective variable. In its full generality the problem is difficult to handle, but the model cases where $f(u)=u^p$ if $N\geq 3$ 
 \bel{I-2}
-\Gd u +u^p -m|\nabla u|^q=0.
\ee
 has already been studied in \cite{BVGHV2}. The results therein shows the rich interaction of the phenomena coming eiter from the interaction of the diffusion and the absorption, the diffusion and the reaction or the reaction and the absorption. A variant of this equation  is the Chipot-Weissler equation
  \bel{I-2+}
-\Gd u -u^p +m|\nabla u|^q=0,
\ee
which presents also interesting interaction between the pairing of the three terms. This has been studied in \cite{BVGHV} and \cite{BVV}. For this last equation the main difficulty is to obtain a priori estimates of positive solutions near an isolated singularity. These estimates have been obtained in 
\cite{BVGHV} and \cite{BVV}.
The aim of the present paper is to study a 2-dimensional problem of type $(\ref{I-2})$ where naturally the power absorption is replaced by an exponential absorption.  Our model equation is 
\bel{I-1}
{\CE}(u):=-\Gd u +ae^{bu} -m|\nabla u|^q=0
\ee
in $B_1\setminus\{0\}$ where $B_1$ is the unit disk of $\BBR^2$, $m,a,b$ are positive parameters and  $q>1$.  For $(\ref{I-2})$ and $(\ref{I-0})$ the a priori estimate near an  isolated singularity is  easy to obtain, the difficulty lies in the full description of singularities and the classifications of singular solutions . Three particular equations are naturally associated to equation $(\ref{I-1})$:\\
the Emden's equation
\bel{I-3}
-\Gd u +ae^{bu}=0,
\ee
the Ricatti equation
\bel{I-4}
-\Gd u -m|\nabla u|^q=0,
\ee
and the eikonal equation
\bel{I-5}
ae^{bu} -m|\nabla u|^q=0.
\ee
The predominance of a specific equation of this type depends on the value of $q$, with critical value $q=2$. When $1<q<2$ the singularities are similar to the one of Emden's equation, an equation which has been completely described in \cite{Va1}.  We recall the main result therein: \\

\nind{\bf V\'azquez theorem }{\it Let $\Omega\subset\BBR^2$ be a bounded smooth domain containing $0$ and  $a,b>0$. Then for any $\phi\in C^1(\prt\Gw)$ and $\gamma\in [0,\frac 2b]$ there exists a unique function $v_\gamma\in C^1(\overline \Gw\setminus\{0\})$ such that $e^{bv_\gamma}\in L^1(\Omega)$ satisfying 
\begin{equation}\label{I-Va}
-\Delta v_\gamma +ae^{bv_\gamma}=2\gp\gamma\gd_0\quad\text{in }{\CD}'(\Omega).
\end{equation}
coinciding with $\phi$ on $\prt\Gw$. Furthermore
\begin{equation}\label{I-Va-1}
\lim_{x\to 0}\frac{v_\gamma(x)}{-\ln |x|}=\gamma.
\end{equation}}
The condition $0\leq\gamma\leq\frac 2b$ is also a necessary condition for the solvability of $(\ref{I-Va})$\\
When $q=2$, the change of unknown $u=\frac 1m\ln v$ transform $(\ref{I-1})$ into 
\bel{I-6}
-\Gd v+mav^{\frac bm+1}=0,
\ee
an equation which has been thoroughly studied (e.g. \cite{BrVe}, \cite{Ve0}) and thus we will omit it. When $q>2$ the dominent equation is the eikonal one. There is a particular explicit solution which plays  a key role in the studies of singularities
\begin{equation}\label{I-7}
U(|x|)=\frac qb\left(\frac ma\right)^{\frac 1q}\ln\left(\frac{q}{b|x|}\right).
\end{equation}
We summarise below our three main results\\

\nind{\bf Theorem 1 }{\it Let $m,a,b>0$ and $1<q<2$. If $u\in C^1(\overline B_1\setminus\{0\})$ is a positive function satisfying $(\ref{I-1})$. Then there exists $\gamma\in [0,\frac 2b]$ such that 
\begin{equation}\label{I-8}
\lim_{x\to 0}\frac{u(x)}{-\ln|x|}=\gamma.
\end{equation}
Furthermore $e^{bu}$ and $|\nabla u|^q$ are integrable in $B_1$ and there holds
\begin{equation}\label{I-9}
-\Delta u +ae^{bu}-m|\nabla u|^q=2\gp\gd_0\quad\text{in }{\CD}'(B_1).
\end{equation}
If $\gamma=0$, $u$ can be extended as a $C^1(\overline B_1)$ function.
}\\

Conversely we prove \\

\nind{\bf Theorem 2 }{\it 
Under the assumptions on $m,a,b,q$ of Theorem 1, for any $\gamma\in (0,\frac 2b]$ and $\phi\in C^1(\prt B_1)$ there exists a unique function $u\in C^1(\overline B_1\setminus\{0\})$  satisfying $(\ref{I-9})$ and coinciding with $\phi$ on $\prt B_1$. Furthermore $e^{bu}$ and $|\nabla u|^q$ are integrable 
in $B_1$ and $(\ref{I-9})$ holds.}\medskip

When $q>2$, the classification of the solutions is totally different and present some similarity with the study of isolated singularities of Lane-Emden equation
\bel{I-10}
-\Gd u-u^p=0,
\ee
in the supercritical case $N\geq 3$ and $p>\frac{N}{N-2}$ (see \cite{GS}). We prove the following characterisation of singular solutions,\\

\nind{\bf Theorem 3 }{\it Let $m,a,b>0$ and $q>2$. If $u\in C^1(\overline B_1\setminus\{0\})$ is a positive solution of $(\ref{I-1})$, then \smallskip

\nind (i) either 
\begin{equation}\label{I-11*}
\lim_{x\to 0}\frac{u(x)}{-\ln|x|}= \frac qb,
\end{equation}

\nind (ii) or $u$ can be extended to $\overline B_1$ as a H\"older continuous function. \smallskip
}\\

\nind {\bf Acknowledgement}: Yimei Li is supported by the National Natural Science Foundation of China (12101038) and the China Scholarship Council (No. 202307090023).
The authors are grateful to M. F. Bidaut-V\'eron for useful discussions during the preparation of this work.
 
\mysection{A priori estimates}
 
 The next estimates are valid in any dimension
 \bprop{A priori} Let $m,a,b>0$ and $q>1$. Assume that for $x_0\in B_1\subset\BBR^N$ and $0<R<|x_0|$ $u$ is a solution of $(\ref{I-1})$ in $B_R(x_0)$. Then there exists 
  $\Lambda=\Lambda (m,q,N,a,b)>0$ such that 
 \smallskip
 
 \nind (i) If $1<q<2$, 
 \begin{equation}\label{I-11}
u(x_0)\leq\frac{2}{b}\ln\frac{1}{|x_0|}+\Lambda.
\end{equation}

\noindent (ii) If $q>2$
  \begin{equation}\label{I-12}
u(x_0)\leq\frac{q}{b}\ln\frac{1}{|x_0|}+\tilde \Lambda.
\end{equation}
 \es
 \Proof Let $R'<R$. We set for $\lambda,\mu>0$, 
$$\psi(x)=\lambda\ln(1/(R'^2-|x-x_0|^2))+\mu=\lambda\ln(1/(R'^2-r^2))+\mu,
$$
with $r=|x-x_0|$. Then 
$$ae^{b\psi}=\frac{ae^{b\mu}}{(R'^2-r^2)^{b\lambda}},
$$
$$|\nabla \psi(x)|=|\psi_r(r)|=\frac{2\lambda r}{R'^2-r^2},
$$
$$|\nabla \psi(x)|^q=\frac{2^q\lambda^q r^q}{(R'^2-r^2)^q},
$$
$$-\Delta \psi(x)=-2\lambda \frac{NR'^2-(N-2)r^2}{(R'^2-r^2)^2}.
$$
Hence there holds in $B_{R'}(x_0)$,
$${\CE}(\psi)=\frac{ae^{b\mu}}{(R'^2-r^2)^{b\lambda}}-2\lambda \frac{NR'^2-(N-2)r^2}{(R'^2-r^2)^2}-\frac{m2^q\lambda^q r^q}{(R'^2-r^2)^q}.
$$
We have $\psi>u$ in a neighbourhood of $\partial B_{R'}(x_0)$. We encounter two cases\smallskip

\noindent (i) $1<q\leq2$. Then we take $b\lambda=2$ and get
\begin{equation}\label {AE-1}\begin{array}{lll}\displaystyle {\CE}(\psi)=(R'^2-r^2)^{-2}\left[ae^{b\mu}-\frac{4}{b}\left(NR'^2-(N-2)r^2\right)-\frac{4^qmr^q}{b^q}\left(R'^2-r^2\right)^{2-q}\right]\\[3mm]
\phantom{\displaystyle {\CE}(\psi)}\displaystyle \geq (R'^2-r^2)^{-2}\left[ae^{b\mu}-\Lambda^* \max\{R'^2,R'^{4-q}\}\right]
\\[3mm]
\phantom{\displaystyle {\CE}(\psi)}\displaystyle \geq (R'^2-r^2)^{-2}\left[ae^{b\gm}-\Gl^*R'^2\right],
\end{array}\end{equation}
where $\Lambda^*=\Lambda^*(N,q,m,b)>0$. We choose
\begin{equation}\label {AE-2}
\mu=\frac{2}{b}\ln R'+\frac{1}{b}\ln\left(\frac{\Lambda^*}{a}\right)
\end{equation}
and derive that ${\CE}(\psi)\geq 0$ in $B_{R'}(x_0)$. Hence $u\leq \psi$ therein. This implies
\begin{equation}\label {AE-3}
u(x_0)\leq \psi(x_0):=\frac 2b\ln\frac 1{R'}+\frac 1b\ln\left(\frac{\Lambda^*}{a}\right).
\end{equation}
Letting $R'\uparrow R=|x_0|$ yields  $(\ref{I-11})$.
\smallskip

\noindent (ii) $q>2$. Then we  take $b\lambda=q$. Hence
\begin{equation}\label {AE-4}\begin{array}{lll}
\displaystyle {\CE }(\psi)=(R^2-r^2)^{-q}\left[ae^{b\mu}-\frac{q}{b}\left(NR^2-(N-2)r^2\right)\left(R^2-r^2\right)^{q-2}-\frac{(2q)^qmr^q}{b^q}\right]\\[3mm]
\phantom{\displaystyle {\CE }(\psi)}\displaystyle \geq (R^2-r^2)^{-q}\left[ae^{b\mu}-\tilde\Lambda^* R^q\right]
\end{array}\end{equation}
where $\tilde\Lambda^*=\tilde\Lambda^*(N,q,m,a,b)>0$. We take
\begin{equation}\label {AE-5}
\mu=\frac{q}{b}\ln R+\frac{1}{b}\ln\left(\frac{\tilde \Lambda^*}{a}\right),
\end{equation}
and derive as in case (i)
\begin{equation}\label {AE-6}
u(x_0)\leq \psi(x_0):=\frac qb\ln\left(\frac 1{R'}\right)+\frac{1}{b}\ln\left(\frac{\tilde \Lambda^*}{a}\right).
\end{equation}
Letting $R'\uparrow R$ we obtain $(\ref{I-12})$.\qeda\medskip

In the next lemma we give some precise estimates on the function $v_\gamma$ solution of $(\ref{I-Va})$
 \blemma {Lem1} Let $m,a,b>0$ and $0<\gamma\leq\frac 2b$. If $v_\gamma\in C^1(\overline B_1\setminus\{0\})$ is the solution $(\ref{I-Va})$ which vanishes on $\prt B_1$, then \smallskip
 
 \nind 1- If $0\leq\gg<\frac 2b$, there exists $\phi^*\in C^1(\overline B_1)$ such that 
 \begin{equation}\label {AE-7*}
v_\gg(x)=\gg\ln\frac{1}{|x|}-\phi^*(x)\quad\text{for all }x\in \overline B_1\setminus\{0\}.
\end{equation}

\nind 2- If $\gg=\frac 2b$, then for all $x\in \overline B_1\setminus\{0\}$ we have\smallskip
 
 \begin{equation}\label{AE-8*}
\frac 2b\left(\ln\frac 1{|x|}-\ln(1-\ln |x|)\right)+\gk_2 \leq v_\frac 2b(x)\leq \frac 2b\left(\ln\frac 1{|x|}-\ln(1-\ln |x|)\right)+\gk_1
\end{equation}
where 
\bel{gk}\gk_1=\max\left\{0,\frac 1b\ln\frac {2}{ab}\right\}\quad\text{and }\;\gk_2=\min\left\{0,\frac 1b\ln\frac {2}{ab}\right\}.
\ee
 In particular 
 \bel{AE9}
 v_\frac 2b(x)= \frac 2b\left(\ln\frac 1{|x|}-\ln(1-\ln |x|)\right),
 \ee
 if $ab=2$.
 \es
 \Proof Let $\gg\in (0,\frac 2b]$ and $\phi_1$ be the first eigenfunction of $-\Delta$ in $W^{12}_0(B_1)$, with corresponding eigenvalue $\lambda_1>0$ and supremum equal to $1$ (at $r=0$). For $\theta>0$ and $A>0$, we set 
\begin{equation}\label{E4*}
h_{A}(x)=\gamma \ln\frac 1{|x|}-A|x|^\theta\phi_1(x)=\gamma \ln\frac 1{r}-Ar^\theta\phi_1(r)\quad\text{with } r=|x|.
\end{equation}
Then there holds in $\CD'(B_1)$, 
$$-\Delta h_A(r)=2\pi\gg\gd_0+A\theta^2r^{\theta-2}\phi_1(r)-A\lambda_1r^{\theta}\phi_1(r)+2A\theta r^{\theta-1}\phi'_1(r),
$$
thus
$$\begin{array}{lll}-\Delta h_A+ae^{bh_A}=2\pi\gg\gd_0+r^{\theta-2}\left(A\theta^2\phi_1(r)-A\lambda_1r^{2}\phi_1(r)+2A\theta r\phi'_1(r)+ar^{2-\theta-b\gamma}e^{-bAr^\theta\phi_1(r)}\right).
\end{array}$$
 {\it Step 1:  the case $b\gamma<2$}\\
 (i) {\it Estimate from above. }
Take $\theta=2-b\gamma$. Since 
$$ae^{-bAr^{2-b\gamma}\phi_1(r)}\geq a-abAr^{2-b\gamma}\phi_1(r),$$ 
we derive that for all $0<|x|<1$,
$$\begin{array}{lll}-\Delta h_A+ae^{bh_A}=r^{-b\gamma}\left(A(2-b\gamma)^2\phi_1(r)-A\lambda_1r^{2}\phi_1(r)+2A(2-b\gamma) r\phi'_1(r)+ae^{-bAr^{2-b\gamma}\phi_1(r)}\right)\\[2mm]
\phantom{-\Delta h_A+ae^{bh_A}}\geq r^{-b\gamma}\left(A(2-b\gamma)^2\phi_1(r)-A\lambda_1r^{2}\phi_1(r)+2A(2-b\gamma) r\phi'_1(r)+a-abAr^{2-b\gamma}\phi_1(r)\right).
\end{array}$$
Then there exists $A_0>0$ such that for $A\leq A_0$ we have 
\begin{equation}\label{E5*}
-\Delta h_A+ae^{bh_A}\geq 2\gamma\pi\delta_0\quad\text{in } {\bf D}'(B_1),
\end{equation}
and $A_0$ can be taken independent of $\gamma$. We derive by the comparison principle that
\begin{equation}\label{E6}
v_\gamma(x)\leq \gamma \ln\frac 1{|x|}-A|x|^{2-b\gamma}\phi_1(x) \quad\text{in } B_1\setminus\{0\}.
\end{equation}
(ii) {\it Estimate from below. } 
We give an explicit expression for $v_\gamma$. We already know that $v_\gamma (x) <\gamma\ln \frac1{|x|}$. If we write
\begin{equation}\label{E*}v_\gamma (x) =\gamma\ln \frac1{|x|}+\phi^*(x),
\end{equation}
then $\phi^*$ vanishes on $\partial B_1$ and satisfies
$$-\Delta\phi^*+ar^{-b\gamma}e^{b\phi^*}=0.
$$
The convex functional 
$$J(\phi)=\int_{B_1}\left(\frac 12|\nabla\phi|^2+\frac a{b|x|^{b\gamma}}e^{b\phi}\right) dx
$$
is well defined on the closed subset $K_-$ of $W^{1,2}_0(B_1)$ of nonpositive functions and it satisfies
$$J(\phi)\geq \frac 12\int_{B_1}|\nabla\phi|^2dx.
$$
Hence it achieves its minimum in $K_-$ at $ \phi^*$. There holds
$$-\Delta (\gamma \ln\tfrac 1{|x|}-A|x|^{2-b\gamma})=A\theta^2|x|^{-b\gamma}
$$
and $\phi^*$ vanishes on $\partial B_1$. Since $|x|^{2-b\gamma}\in L^{1+\epsilon}(B_1)$, $\phi^*\in W^{2,1+\epsilon}( B_1)\subset C^{0,\frac{2\ge}{1+\ge}}(\overline B_1)$.
Consequently
\begin{equation}\label{E7}
v_\gamma(x)= \gamma \ln\frac 1{|x|}+\phi^*(x) \quad\text{in } B_1\setminus\{0\},
\end{equation}
and
\begin{equation}\label{E7*}
v_\gamma(x)-\gamma \ln\frac 1{|x|}\to \phi^*(0)<0 \quad\text{as } x\to 0.
\end{equation}

 \nind {\it Step 2:  the case $b\gamma=2$.} For $\gk\in\BBR$ we set
$$ \psi_{\gk}(r)=\frac 2b\left(\ln \frac 1r-\ln(1-\ln r)\right)+\gk.
$$
Then $\psi_{\gk}(1)=\gk$ and 
\begin{equation}\label{N1}-\Delta \psi_{\gk}+ae^{b\psi_{\gk}}=\frac{1}{r^2}\left(-\frac{2}{b(1-\ln r)^2}+\frac{ae^{b\gk}}{(1-\ln r)^2}\right).
\end{equation}
Since 
$$\int_{B_1}\frac{dx}{|x|^2(1-\ln |x|)^2}=2\gp \int_0^1\frac{dr}{r(1-\ln r)^2}=2\gp$$
the function $e^{b\psi_\gk}$ is integrable in $B_1$ and $\psi_{\gk}$ satisfies in $\CD'(B_1)$
\begin{equation}\label{N2}
-\Delta\psi_\gk+ae^{b\psi_\gk}=\left(ae^{b\gk}-\frac{2}{b}\right)\frac{1}{|x|^2(1-\ln |x|)^2}+\frac{4\pi}{b}\delta _0.
\end{equation}
Hence $\psi_\gk$ is a supersolution nonnegative on $\prt B_1$ if we choose $\gk=\gk_1$ where
\begin{equation}\label{kappa1}
\gk_1=\max\left\{0, \frac{1}{b}\ln\frac{2}{ab}\right\}
\end{equation}
It is a subsolution if we take $\gk=\gk_2$ where
\begin{equation}\label{kappa2}
\gk_2=\min\left\{0, \frac{1}{b}\ln\frac{2}{ab}\right\}.
\end{equation}
Therefore
\begin{equation}\label{N'2} \frac 2b\left(\ln \frac 1{|x|}-\ln(1-\ln |x|)\right)+\gk_2\leq v_{\frac 2b}(x)\leq \frac 2b\left(\ln \frac 1{|x|}-\ln(1-\ln |x|)\right)+\gk_1,
\end{equation}
with equality if $b=\frac 2a$ since in that case $\gk_1=\gk_2=0$.

\qeda
\mysection{Singular solutions in the case $1<q<2$}
\bdef{singdir} A function $u\in C^1(\overline B_1\setminus\{0\})$ is a solution of the singular Dirichlet problem
\begin{equation}\label{Ex1}\BA{lll}
-\Gd u+ae^{bu}-m|\nabla u|^q=2\pi\gg\delta_0\qquad&\text{in }\CD'(B_1)\\
\phantom{-\Gd +ae^{bu}-m|\nabla u|^q}
u=\phi\qquad&\text{on }\prt B_1,
\EA\end{equation}
where $\phi$ is a $C^1$ function, if $u$, $e^{bu}$ and $|\nabla u|^q$ are integrable in $B_1$ and 
\begin{equation}\label{W1}
\int_{B_1}\left(-u\Gd\gz+\left(ae^{bu}-m|\nabla u|^q\right)\gz\right)dx=2\pi\gg\gz(0)\quad\text{for all }\,\gz\in C^2_c(B_1).
\end{equation}
\es
In view of Vazquez theorem, any solution $u_\gg$ of $(\ref{Ex1})$ is bounded from below by the solution $v_\gg$

\subsection{Radial solutions}
 The expression of $(\ref{I-2})$ for a radial function $u:=u(r)$ defined in $B_1\setminus\{0\}$ is 
 \begin{equation}\label{Rad1}
 -u''-\frac 1ru'+ae^{bu}-m|u'|^q=0\quad\text{for }0<r<1.
 \end{equation}
As a first gradient estimate for radial solutions we prove
\blemma{grad1} Let $1<q<2$ and $u$ be a nonnegative  solution of $(\ref{Rad1})$  in $B_1\setminus\{0\}$. Then for any $0<r_0<1$
 there exists $c=c(m,q)>0$ such that 
 \begin{equation}\label{E23-4}
|u_r(r)|\leq \frac cr\qquad\text{for all }0<r\leq r_0.
\end{equation}
\es
\Proof Clearly $u$ cannot have a local maximum and up to a translation (it changes only the value of $a$) we can assume that $u$ vanishes on $\partial B_1$, hence it is decreasing. There holds
 $$-u''-\frac 1ru'-m|u'|^q\leq 0,
 $$
 then
 $$\left(\frac 1{1-q}|u'|^{1-q}\right)'+\frac{1}{r}|u'|^{1-q}\leq m,
 $$
 and thus
 $$\left(r^{1-q}|u'|^{1-q}+m\frac{q-1}{2-q}r^{2-q}\right)'\geq 0.
 $$
 Therefore the function $r\mapsto r^{1-q}|u'|^{1-q}+m\frac{q-1}{2-q}r^{2-q}$ is increasing and nonnegative, and thus it admits a limit $\ell\geq 0$  when 
 $r\to 0$. This limit is the same as the one of $r\mapsto r^{1-q}|u'|^{1-q}$. If $\ell>0$ 
 $$r^{1-q}|u'|^{1-q}= \ell (1-o(1))\quad\text{when }\,r\to 0.
 $$
Hence 
$$\lim_{r\to 0}r^{1-q}|u'|^{1-q}= \ell
$$
which implies $(\ref{E23-4})$. If $\ell=0$, then 
$$\lim_{r\to 0}r|u'(r)|= \infty.
$$
This implies that for all $A>0$, one has
$$\lim\inf_{r\to 0}\left(u(r)-A\ln \frac{1}{r}\right)=\infty.
$$
This contradicts the a priori estimate obtained before. \qeda\medskip

The next easy to prove result will be used several times in the sequel.
\blemma{estimL} Let $\gg\geq 0$, $F\in L^1(B_1)$, $\phi\in C(\prt B_1)$ and $w\in L^1(B_1)\cap C(\overline B_1\setminus\{0\})$ be the weak solution of
\bel{ELL1}\BA{lll}
-\Gd w=F+2\gp\gg\gd_0\quad&\text{in }{\CD}'(B_1)\\
\phantom{-\Gd }
w=\phi\quad&\text{in }\prt B_1.
\EA\ee
Then there holds
\bel{ELL2}
\bar w(r)=\gg\ln\frac 1r+o(\ln r)\quad\text{as }|x|=r\to 0
\ee
where $\dsps\bar w(r)=\frac 1{2\gp}\int_{0}^{2\gp}w(r,\gth)d\gth$. If we assume that $F\in L^{1+\gd}(B_1)$ for some $\gd>0$, then 
\bel{ELL2+}
w(x)=\gg\ln\frac 1{|x|}+O(1)\quad\text{as }|x|\to 0
\ee
\es
\Proof We write $w(x)=\gg E(x)+\psi(x)$ where $E(x)=\ln\frac 1{|x|}$ and $\psi$ is the solution of 
\bel{ELL3}\BA{lll}
-\Gd \psi=F\quad&\text{in }{\CD}'(B_1)\\
\phantom{-\Gd }
\psi=\phi\quad&\text{in }\prt B_1.
\EA\ee
Let  $\bar \psi(r)$, $\bar \phi$ and $\bar F(r)$ be the circular averages of $\psi(r,.)$, $\phi(.)$ and $\bar F(r,.)$. Since
$$-(r\bar \psi'(r))'=r\bar F(r),
$$
A standard approximation argument associated to uniqueness yields 
$$\bar\psi(r)=\bar\phi+\int_r^1s^{-1}\int_0^s t\bar F(t)dt ds.
$$
By l'Hospital's rule we obtain $(\ref{ELL2})$. If $F\in L^{1+\gd}(B_1)$, then $\psi\in W^{2,1+\gd}_{loc}(B_1)$. By Morrey-Sobolev imbedding inequality, 
$W^{2,1+\gd}_{loc}(B_1)\subset C^{0,\frac{2\gd}{1+\gd}}(B_1)$. Hence $\psi$ is H\"older continuous in $B_1$, which  ends the proof. \qeda

\blemma{Lem2}
Let  $1<q<2$, $m,a,b>0$ and $\gamma\in \left(0,\frac{2}{b}\right]$. If $u:=u_\gg$ is a radial solution of 
\[
-\Delta u+ae^{bu}-m|\nabla u|^{q}=2\pi\gamma\delta_{0}\quad \mbox{in}\  D'(B_{1})
\]
which vanishes on $\partial B_1$, then \smallskip

\nind 1-  If $0<\gamma<\frac{2}{b}$, there exist a  positive constant $A$ and a  continuous positive function $\phi^{*}$  in $B_{1}$ which vanishes on $\partial B_{1}$ such that in $B_{1}\setminus \{0\}$,
\[
\gamma\ln\frac{1}{|x|}-\phi^{*}(x)\leq u_{\gamma}\leq  \gamma\ln\frac{1}{|x|}-A|x|^{2-b\gamma}\phi_{1}(x)+mc^q\eta(x),
\]
where $\phi_{1}$ is the first eigenfunction of $-\Delta $ in $W^{1,2}_{0}(B_{1})$ and  $\eta$ satisfies

\bel{eta}\BA{lll}
-\Delta \eta=  |x|^{-q}\quad&\mbox{\rm in }\  B_{1}\\
\phantom{-\Delta }\eta=0\quad&\mbox{\rm on }\  \partial B_{1},
\EA
\ee
with $c$ is the constant in $(\ref{E23-4})$, and actually, $\eta(x)=\frac{1}{(2-q)^{2}}(1-|x|^{2-q}).$

\nind 2- If $\gamma=\frac{2}{b}$, we have in $B_{1}\setminus \{0\}$,
\[
\frac{2}{b}\left(\ln\frac{2}{|x|}-\ln(1-\ln|x| )\right)+\gk_2\leq u_{\frac{2}{b}}(x)\leq  \frac{2}{b}\left(\ln\frac{2}{|x|}-\ln(1-\ln|x| )\right)+\gk_1+\eta(x),
\]
where the $\gk_i$ have been defined in $(\ref{gk})$.
\es
\Proof
Because of $(\ref{E23-4})$ $u_{\gamma}$ satisfies
$$
-\Delta u_{\gamma}+ae^{bu_{\gamma}}\leq  mc^{q}|x|^{-q}+2\pi\gamma\delta_{0}\quad\mbox{\rm in }\  D'(B_{1}).
$$
Since $v_{\gamma}+mc^q\eta$ vanishes on $\partial B_1$ and satisfies
$$-\Delta (v_{\gamma}+mc^q\eta)+ae^{b(v_{\gamma}+mc^q\eta)}\geq  mc^{q}|x|^{-q}+2\pi\gamma\delta_{0}\quad\mbox{\rm in }\  D'(B_{1}),
$$
we deduce by the comparison  principle as in \cite{Va1}, that
\bel{comp1}
 u_{\gamma}\leq v_{\gamma}+mc^q\eta\quad\mbox{in} \quad B_{1}\setminus \{0\}.
\ee
On the other hand,  
$$-\Delta u_{\gamma}+ae^{bu_{\gamma}}\geq 2\pi\gamma\delta_{0}.
$$
Since $ae^{bu_{\gamma}}\in L^1(B_1)$, it follows by the comparison  principle again, that
\bel{comp2}
 v_{\gamma}\leq u_{\gamma}\quad\mbox{in} \quad B_{1}\setminus \{0\}.
\ee
Therefore, we conclude that
\bel{comp3}
 v_{\gamma}\leq u_{\gamma}\leq v_{\gamma}+mc^q\eta\quad\mbox{in} \quad B_{1}\setminus \{0\}.
\ee
Combining with \rlemma{Lem1}, we finish the proof.
\qeda

\subsection{General solutions}
In this section we deal with the questions of existence, classification and uniqueness of solutions of $(\ref{Ex1})$.
\subsubsection{Existence of general solutions}
\bth{Exis} Let $1<q<2$ and $m,a,b>0$. Then for any $\gg\in (0,\frac 2b]$ and $\phi\in C^1(\prt B_1)$ there exists a function $u:=u_\gg\in C^1(\overline B_1\setminus\{0\})$ such that 
$e^{bu}$ and $|\nabla u|^q$ are integrable in $B_1$ satisfying $(\ref{Ex1})$ and 
\begin{equation}\label{Ex2*}
\lim_{x\to 0}\frac{u(x)}{-\ln |x|}=\gg
\end{equation}
\es


\nind {\it Proof  of \rth{Exis}.} {\it Step 1: existence.} For the sake of simplicity, we assume $\phi\geq 0$; this can be achieved up to replacing $a$ by some $\tilde a$. The function $v_\gg$ which satisfies $(\ref{I-Va})$ with boundary data $\phi$ is a subsolution in $B_1\setminus\{0\}$. We look for a supersolution under the form
$$\psi(x)=\gg\ln\frac 1{|x|}+c(1-|x|^{2-q})+M,
$$
where $c$ and $M$ are to be found. Since $-\Delta (1-|x|^{2-q})=(2-q)^2|x|^{-q}$
$$\CE(\psi)=a|x|^{-b\gg}e^{ bc(1-|x|^{2-q})+bM}+|x|^{-q}\left(c(2-q)^2-m\gg^q\left(1+\frac{c(2-q)}{\gg}r^{2-q}\right)^q\right)
$$
There exists $r_0\in (0,1)$ and $c_0>0$ such that 
 $$c_0(2-q)^2-m\gg^q\left(1+\frac{c_0(2-q)}{\gg}r^{2-q}\right)^q\geq 0\quad\text{for all }0<r\leq r_0.
 $$
 Then we choose $M>0$ large enough such that 
 $$ae^{bM}\geq r_0^{-q}m\gg^q\left(1+\frac{c_0(2-q)}{\gg}\right)^q,
 $$
 and
 $$\psi(x)=M\geq \norm\gf_{L^\infty(\prt B_1)}.
 $$
 Hence $\CE(\psi)\geq 0$ in $B_1\setminus\{0\}$, and $\psi\geq\phi$ on $\prt B_1$. It follows that $\psi$ is a supersolution, larger than $v_\gamma$. Therefore there exists $u=u_\gamma$ coinciding with $v_\gamma$, such that 
\bel{AA}v_\gamma\leq u_\gamma\leq \psi
\ee
 which satisfies  $(\ref{I-1})$ in $B_1\setminus\{0\}$.\smallskip
 
 \nind {\it Step 2: upper estimate.} The function $u_\gamma$ may not be radial since it is not assumed that $v_\gamma$ is so. But replacing $v_\gamma$ by 
 $v_{rad,\gamma}$ which satisfies the same equation $(\ref{I-Va})$ with boundary data $\norm \phi_{L^\infty(\prt B_1)}$, we have that $v_{rad,\gamma}\geq v_\gamma$
  with boundary data $\norm\gf_{L^\infty(\prt B_1)}$. Hence there exists a radial solution $u_{rad,\gamma}$ of $(\ref{I-1})$ satisfying 
\bel{XX1}v_\gamma\leq v_{rad,\gamma}\leq u_{rad,\gamma}\leq\psi
\ee
 and $u_{rad,\gamma}\geq  u_\gamma$ by the standard construction iterative scheme. By \rlemma {grad1} there holds
 $$|u'_{rad,\gamma}(r)|\leq \frac cr.
 $$
 Therefore
 $$
-\Gd u_{rad,\gamma}+ae^{bu_{rad,\gamma}}=m|\nabla u_{rad,\gamma}|^q\leq mc^q r^{-q},
 $$
 and $(\ref{XX1})$ gives that
 $$\lim_{r\to 0}\frac{u_{rad,\gamma}(r)}{-\ln r}=\lim_{x\to 0}\frac{u_\gamma(x)}{-\ln |x|}=\gamma.
 $$
 Let $\{\theta_n\}\subset C^{1,1}(\BBR^2$) defined by
 $$\theta_n(x)=\left\{\BA{lll}\dsps
 \;2n^2\left(|x|-\frac{1}{n}\right)_+^2\qquad&\text{if }|x|\leq \frac{3}{2n}\\[4mm]\dsps
 -2n^2|x|^2+8n|x|-7\qquad&\text{if }\frac{3}{2n}\leq |x|\leq \frac{2}{n}\\[2mm]
\; 1\qquad&\text{if }|x|\geq \frac{2}{n},
 \EA\right.
 $$
then $0\leq\theta_n\leq 1$, $\theta_n=0$ on $B_{\frac 1n}$, $\theta_n=1$ on $B^c_{\frac 2n}$, $|\nabla\theta_n|\leq cn{\bf 1}_{B_{\frac 2n}\setminus B_{\frac 1n}}$ and
 $$\Gd\theta_n(x)=4n^2\left(2-\frac{1}{n|x|}\right){\bf 1}_{B_{\frac 3{ 2n}}\setminus B_{\frac 1n}}
 -4n^2\left(2-\frac{2}{n|x|}\right){\bf 1}_{B_{\frac 2n}}\setminus B_{\frac 3{ 2n}}.
 $$
 If $\gz\in C_c^{\infty}(B_1)$ is nonnegative with $\gz(0)=1=\max_{|\gz|\leq 1}\gz(z)$.  We have
\bel{XX0}-\int_{B_1}u_{rad,\gamma}\theta_n\Gd\gz  dx+a\int_{B_1}e^{bu_{rad,\gamma}}\gz \theta_n dx-\int_{B_1}\left(\gz\Gd \theta_n+2\nabla \theta_n.\nabla\gz\right) u_{rad,\gamma} dx
 \leq mc^q\int_{B_1}\theta_n\gz|x|^{-q}dx.
 \ee
 Clearly there holds when $n\to\infty$,
 $$\int_{B_1}u_{rad,\gamma}\theta_n\Gd\gz  dx\to \int_{B_1}u_{rad,\gamma}\Gd\gz  dx,
 $$
 and
 $$mc^q\int_{B_1}\theta_n\gz|x|^{-q}dx\to mc^q\int_{B_1}\gz|x|^{-q}dx.
 $$
 Next
 $$\int_{B_1}\left|\nabla \theta_n.\nabla\gz\right| u_{rad,\gamma} dx\leq cn\int_{B_{\frac 2{ n}}\setminus B_{\frac 1n}}u_{rad,\gamma}|\nabla \gz|dx\to 0,
 $$
 since $ u_{rad,\gg}(x)= \gg\ln  (\frac 1{|x|})(1+o(1))$. The last term is more delicate to estimate:
 $$-\int_{B_1}\Gd \theta_n u_{rad,\gamma}\gz dx=-\int_{B_{\frac 3{ 2n}}\setminus B_{\frac 1n}}\Gd \theta_n u_{rad,\gamma}\gz dx-\int_{B_{\frac 2n}\setminus B_{\frac 3{ 2n}}}\Gd \theta_n u_{rad,\gamma}\gz dx:=I+II.
 $$
Since $\gz(x)=1+O(|x|)$ and  $ u_{rad,\gg}$ satisfies $(\ref{XX1})$, estimates on $v_\gg$ and $\psi$ yield
\bel{YY1}\gg\ln\frac1{|x|}-M\leq u_{rad,\gamma}(x)\leq \gg\ln\frac1{|x|}+M
\ee
for some $M>0$, and therefore
$$I=-4n^2\int_{B_{\frac 3{ 2n}}\setminus B_{\frac 1n}}\left(2-\frac{1}{n|x|}\right)u_{rad,\gamma} dx+O\left(\frac{\ln n}{n}\right).
$$
Hence
$$-6\gp\max\left\{u_{rad,\gamma}: \frac{1}{n}\leq |x|\leq \frac{3}{2n}\right\}+O\left(\frac{\ln n}{n}\right)\leq I\leq -6\gp\min\left\{u_{rad,\gamma}: \frac{1}{n}\leq |x|\leq \frac{3}{2n}\right\}+O\left(\frac{\ln n}{n}\right).
$$
In the same way
$$II=4n^2\gz(0)\int_{B_{\frac 2{ n}}\setminus B_{\frac 3{2n}}}\left(2-\frac{2}{n|x|}\right)u_{rad,\gamma}\gz dx+O\left(\frac{\ln n}{n}\right).
$$
Hence
$$6\gp\min\left\{\tilde u_\gg: \frac{3}{2n}\leq |x|\leq \frac{2}{n}\right\}+O\left(\frac{\ln n}{n}\right)\leq II\leq 6\gp\max\left\{\tilde u_\gg: \frac{3}{2n}\leq |x|\leq \frac{2}{n}\right\}+O\left(\frac{\ln n}{n}\right).
$$
Using $(\ref{YY1})$ we see that up to some constants $M>0$
Therefore on
$$-6\gp\ln n-M+O\left(\frac{\ln n}{n}\right)\leq I\leq -6\gp\ln n+M+O\left(\frac{\ln n}{n}\right),
$$
and 
$$6\gp\ln n-M+O\left(\frac{\ln n}{n}\right)\leq II\leq 6\gp\ln n+M+O\left(\frac{\ln n}{n}\right).
$$
Thus $I+II$ remains bounded independently of $n$. Combining these estimates with $(\ref{XX0})$ we conclude that $\int_{B_1}e^{bu_{rad,\gamma}}\gz \theta_n dx$ remains bounded independently of $n$ too. By Fatou's lemma $e^{bu_{rad,\gamma}}\gz\in L^1(B_1)$. Since $u_\gamma\leq u_{rad,\gamma}$ we infer that $e^{b u_\gg}\in L^1(B_1)$. Because 
$$-\Gd u_\gg=m|\nabla u_\gg|^q-ae^{bu_\gg},
$$
it follows by Brezis-Lions lemma that $|\nabla u|^q\in L^1(B_1)$ and there exist $\tilde \gamma\geq 0$ such that 
$$-\Gd u_\gg=m|\nabla u_\gg|^q-ae^{bu_\gg}+2\pi \tilde \gamma\delta_0\quad\text{in }\CD'(B_1).
$$
Set $F=m|\nabla u_\gg|^q-ae^{bu_\gg}$, Since $F\in L^1(B_1)$ it follows by \rlemma{estimL} that
$$\bar u_\gg(r)=-\tilde \gg\ln{r} +o(\ln r).
$$
Combined with $(\ref{AA})$ we obtain that $\tilde\gg=\gg$.\qeda
\medskip

The next estimate compares $v_\gg$ and $u_\gg$ 

\bprop{Uniq+exis1} Let $\gg\in (0,\frac 2b]$ and $\gf\in C^1(\partial B_1)$ is nonnegative. If  $u\in C^1(\overline B_1\setminus\{0\})$ satisfies $e^{bu}\in L^1(B_1)$ and $|\nabla u|^q\in L^1(B_1)$ is a solution of $(\ref{Ex1})$, then there exists $K>0$ such that 
\bel{sp9}
v_\gg\leq u(x)\leq v_\gg+K\quad\text{for }x\in B_1\setminus\{0\},
\ee
where $v_\gg$ is the solution of $(\ref{I-Va})$ which vanishes on $\prt B_1$.
\es
\Proof Since $u$ is a solution, it is a supersolution of $(\ref{I-Va})$, and by the order principle it is larger than solution $v_\gg$ with the same boundary value. Using the regularity estimates of the Laplacian in the space of bounded Radon measures $\frak M_b(B_1)$ we have that $\nabla u\in M^2(B_1)$ where $M^2(B_1)$ is the Marcinkiewicz space 
coinciding with the Lorentz space $L^{2,\infty}(B_1)$. What is important to know is that $M^2(B_1)\subset L^{2-\ge}(B_1)$ for any $\ge>0$. Therefore 
$|\nabla u|^q\in L^{\frac{2-\ge}{q}}(B_1)$. Since $q<2$ we take $\ge>0$ such that $\frac{2-\ge}{q}=1+\gd>1$. Let $u_1$ be the solution of 
$$\BA{lll}-\Gd u_1=m|\nabla u|^q&\quad\text{in } B_1\\
\phantom{-\Gd}
u_1=\gf&\quad\text{on } B_1.
\EA
$$
The function $u_1$ is nonnegative and belongs to $W^{2,1+\gd}(B_1)$. By Morrey-Sobolev imbedding theorem, $W^{2,1+\gd(B_1)}\subset C^{0,\frac{2\gd}{1+\gd}}(\overline B_1)$. This implies that $u_1$ is bounded from above by some $K>0$. The function $v_\gg+u_1$ satisfies 
$$-\Gd(u_1+v_\gg)+ae^{b(u_1+v_\gg)}>-\Gd(u_1+v_\gg)+ae^{bv_\gg}= m|\nabla u|^q+2\gp\gg\gd_0.
$$
Since $ae^{b(u_1+v_\gg)}\in L^1(B_1)$ we can apply Vazquez' monotonicity estimate \cite{Va1} between $u_1+v_\gg$ and $u$ since the right-hand side of the respective equations coincide and belong to  $\frak M_b(B_1)$. Therefore 
\bel{sp10}
u(x)\leq v_\gg(x)+K,
\ee
and $(\ref{sp9})$ follows.
\qeda \medskip

As a consequence of \rlemma{Lem1} we obtain global estimates on the solutions

\bcor{glob} Under the assumptions of \rprop{Uniq+exis1} any function  $u\in C^1(\overline B_1\setminus\{0\})$ such that $e^{bu}\in L^1(B_1)$ and $|\nabla u|^q\in L^1(B_1)$ of solution of $(\ref{Ex1})$ satisfies, for some $K\geq 0$, the following two-sided inequalities for all $x\in \overline B_1\setminus\{0\}$,\smallskip

 \nind 1- If $0\leq\gg<\frac 2b$, then
 \begin{equation}\label {AE-7}
\gg\ln\frac{1}{|x|}-K\leq u(x)\leq \gg\ln\frac{1}{|x|}+K.
\end{equation}

\nind 2- If $\gg=\frac 2b$, then 
 \begin{equation}\label{AE-8}
\frac 2b\left(\ln\frac 1{|x|}-\ln(1-\ln |x|) \right)-K\leq u(x)\leq \frac 2b\left(\ln\frac 1{|x|}-\ln(1-\ln |x|) \right)+K.
\end{equation}
\es





%



\subsubsection{Classification of isolated singularities}
In this section we prove the following
\bth{class} Let $1<q<2$, $m,a,b>0$. If $u\in C^1(\overline B_1\setminus\{0\})$ is a solution of $(\ref{I-1})$ in $B_1\setminus\{0\}$, then 
there exists $\gg\in [0,\frac{2}{b}]$  such that
\bel{lim1}
\lim_{x\to 0}\frac{u(x)}{-\ln |x|}=\gg.
\ee
Furthermore $e^{bu}\in L^1(B_1)$, $|\nabla u|^q\in L^1(B_1)$ and $u$ is a solution of $(\ref{Ex1})$.\es
\Proof We can assume that $u\lfloor_{\prt B_1}=\phi\in C^{1}(\prt B_1)$ is nonnegative. \\
 {\it 1}- If $e^{bu}\in L^1(B_1)$, it follows by Brezis-Lions Lemma \cite {BrLi} that $|\nabla u|^q\in L^1(B_1)$ and there exists $\gg\geq 0$ such that $(\ref{Ex1})$ holds. By \rlemma {estimL} the circular average $\bar u$ of $u$ satisfies
 $$\bar u(r)=\gg\ln\frac 1r+o(\ln r)\quad\text{as }r\to 0.
 $$
 By convexity
 $$\overline {e^{bu}}(r)\geq e^{b\bar u(r)}.
 $$
Since for any $\ge>0$ there holds
$$\bar u(r)\geq (\gg-\ge)\ln\frac{1}{r}-K_\ge,
$$  
we obtain
$$e^{b\bar u(r)}\geq \frac{ e^{-bK_\ge}}{|x|^{b(\gg-\ge)}}.
$$
If $\gg>\frac 2b$ we can choose $\ge>0$ such that $b(\gg-\ge)>2$. Hence $e^{b\bar u(.)}$ is not integrable in $B_1$, neither $e^{bu}$, contradiction. This implies that 
$\gg\in [0,\frac 2b]$. By \rprop{Uniq+exis1} and  \rcor{glob} we obtain $(\ref{lim1})$.


 
\nind {\it 2}- Next we assume that $e^{bu}\notin L^1(B_1)$, then using the inequality satisfied by $\bar u$,
\bel{ZZ1}r\bar u'(r)+a\int_r^1\overline {e^{bu}}tdt= \bar u'(1)+m\int_r^1\overline {|\nabla u|^q}tdt
\ee
where $\bar u$, $\overline{E^{bu}}$ and $\overline {|\nabla u|^q}$ denote the averages of  $ u$, ${E^{bu}}$ and $ {|\nabla u|^q}$ on $S^1$. Next we will distinguish according $|\nabla u|^q$ is or is not integrable in $B_1$, although if $u$ is radial \rlemma {grad1} implies that $|\nabla u|^q\in L^1(B_1)$.

\nind {\it 2-(i)}- If $|\nabla u|^q\in L^1(B_1)$ we obtain from $(\ref{ZZ1})$ that 
$$\lim_{r\to 0}r\bar u'(r)=-\infty.
$$
By integration, it implies that 
$$\lim_{r\to 0}\frac{\bar u(r)}{-\ln r}=\infty,
$$
which again contradicts the a priori estimate $(\ref{I-11})$. 

\nind {\it 2-(ii)}- If $|\nabla u|^q\notin L^1(B_1)$, then for any $\ge\in (0,1)$,
$$\int_{B_\ge}|\nabla u|^qdx=\infty.
$$
Then for any $k>0$ there exists $\gt=\gt(\ge,k)$ such that 
$$m \int_{B_\ge}|\nabla \min\{u,\gt\}|^qdx=2\gp k.
$$
Since 
$$-\Gd\min\{u,\gt\}+ae^{b\min\{u,\gt\}}\geq m|\nabla \min\{u,\gt\}|^q{\bf 1}_{B_\ge}\quad \text {in }B_1
$$
$\min\{u,\gt\}$ is bounded from below by the solution $v_{k,\ge}$ of 
$$-\Gd v_{k,\ge}+ae^{b v_{k,\ge}}= m|\nabla \min\{u,\gt\}|^q{\bf 1}_{B_\ge}\quad \text {in }B_1.
$$
We have $m|\nabla \min\{u,\gt\}|^q{\bf 1}_{B_\ge}\to 2\gp k\gd_0$ in the weak sense of bounded Radon measures in $B_1$.
By Vazquez' stability theorem, if $k\leq \frac{2}{b}$, $v_{k,\ge}\to v_k$ as $\ge\to 0$. Taking in particular $k=\frac 2b$ we obtain
\bel{Z1}u(x)\geq v_{\frac 2b}(x)\quad\text{for all } x\in B_1\setminus\{0\}.
\ee
Combining this estimate with $(\ref{I-11})$ and the estimate on $v_{\frac 2b}$ yields, for some $K\geq 0$,
\bel{Z2}\frac 2b\left(\ln\frac 1{|x|}-\ln(1-\ln |x|)\right)-K\leq u(x)\leq \frac 2b\ln\frac 1{|x|}+K\quad\text{for all } x\in B_1\setminus\{0\}.
\ee
We set $w(t,.)=u(r,.)+\frac 2b\ln r$ with $t=\ln r$. Then $w$ satisfies in $(-\infty,0)\ti S^1$
\bel{EQ}w_{tt}+w_{\theta\theta}-ae^{bw}+me^{(2-q)t}\left(\left(\frac 2b-w_t\right)^2+w_\gth^2\right)^\frac q2=0,
\ee
and the two-sided estimate
\bel{EQ-1}
-\frac{2}{b}\ln(1-t)-K\leq w(t,.)\leq K.
\ee
Thus the function  $z(t,.)=\frac{w(t,.)}{\ln (1-t)}$
satisfies in $(-\infty,0)\ti S^1$,
\bel{sp3}\BA{lll}\dsps
z_{tt}+z_{\gth\gth}-\frac{2}{(1-t)\ln(1-t)}z_t-\frac{1}{(1-t)^2\ln(1-t)}z-\frac{a}{\ln (1-t)}e^{b\ln(1-t)z} \\[4mm]\phantom{---}\dsps
+me^{(2-q)t}\ln^{q-1} (1-t)\left(\left(\frac 2{b\ln (1-t)}-z_t+\frac z{(1-t)\ln (1-t)}\right)^2+ z_\gth^2\right)^{\frac q2}=0,
\EA\ee
and ($\ref{EQ-1}$) yields for some $T<0$,
$$-\frac 2b+o(1)\leq z(t,.)\leq \frac{K}{\ln(1-t)}\leq o(1)\quad\text{in }(-\infty,T]\ti S^1.
$$
Since $z$ is bounded and negative and all the coefficients are smooth and bounded, the standard regularity theory implies that 
$z$ is bounded in $C^{2,\ga}([T'-1,T'+1]\ti S^1)$. Since 
$w(t,.)=z(t,.)\ln (1-t)$ we have also
\bel{sp3+}
\left(w^2_t(t,\gth)+w^2_\gth(t,\gth)\right)^{\frac 12}\leq c\ln(1-t).
\ee
Returning to $u(r,\gth)$ implies 
\bel{Z3}|\nabla u(x)|=|\nabla u(r,\gth)|=(u_r^2+r^{-2}u_\gth^2)^{\frac 12}\leq c\frac{\ln(1-\ln|x|)}{|x|}\quad\text{for all } x\in B_1\setminus\{0\}.
\ee
Since $q<2$, the right-hand side of $(\ref{Z3})$ belongs to $L^q(B_1)$ and we get a contradiction. This ends the proof.
\qeda

\subsubsection{Uniqueness}
Sharp estimates of the gradient of singular solutions near $x=0$ are  the key-step for proving their uniqueness. 


\bth{Uniq-1} Let $q,m,a,b$ be as in \rth{Exis}. Then for any $\gg\in (0,\frac 2b]$ and $\phi\in C^1(\partial B_1)$ there exists a unique function 
$u\in C^1(\overline B_1\setminus\{0\})$ satisfying equation $(\ref{I-1})$ in $B_1\setminus\{0\}$ and
\bel{Ex2}
u(x)=\gg\ln\frac{1}{|x|}\left(1+o(1)\right).
\ee
Furthermore 
$e^{bu}$ and $|\nabla u|^q$ are integrable in $B_1$, and $u$ satisfies $(\ref{Ex1})$. 
\es
\Proof Let $u$ be any solution nonnegative solution of $(\ref{I-1})$ in $B_1\setminus\{0\}$ satisfying $(\ref{Ex2})$ (the sign requirement on $u$ is not a restriction since we can change $a$). By \rth{class} $e^{bu}$ and $|\nabla u|^q$ are integrable in $B_1$ and $u$ satisfies $(\ref{Ex1})$. Therefore the estimates of \rcor{glob} hold.
\smallskip

\nind  {\it 1- The case $0<\gg<\frac 2b$}.\smallskip

\nind {\it 1-(i) Convergence and gradient estimates.} Under the assumption $(\ref{Ex2})$ we first show that there holds
 \begin{equation}\label{EX2**}
\lim_{|x|\to 0}\left(u(x)-\gg\ln\frac{1}{|x|}\right)=\ell\quad\text{and }\,  \lim_{|x|\to 0}|x|\nabla u(x)=-\gg {\bf i}
\end{equation}
for some $\ell\in [-K,K]$ and ${\bf i}=\frac x{|x|}$. For the proof, we set 
$$w(t,\gth)=u(r,\gth)-\gamma\ln \frac 1r,
$$
with $t=\ln r\in (-\infty,0]$ and $\gth\in S^1$. Then $w$ is bounded in $(-\infty,0]\ti S^1$ where there holds

\bel{EQ*}w_{tt}+w_{\theta\theta}-ae^{(2-b\gg)t}e^{bw}+me^{(2-q)t}\left(\left(\gg-w_t\right)^2+w_\gth^2\right)^\frac q2=0.
\ee
Since $2-b\gg> 0$ and  $2-q>0$, by the standard regularity theory \cite [Theorem A1]{LaLi}, we have the estimate 
$$\BA{lll}\max\{|w_t(t,\theta)|+|w_\gth(t,\theta)|: (t,\theta)\in [T-1,T+1]\ti S^1 \}\\[2mm]
\phantom{\max\{: (t,\theta)\in [T-1,T+1]\ti S^1 \}}\leq c\max\{|w(t,\theta)|: (t,\theta)\in [T-2,T+2]\ti S^1 \}\leq cK,
\EA$$
for all $T\leq -2$. Equivalently $|u_\gth|+|ru_r+\gg|$ is bounded. Now, by linear elliptic equation regularity theory $w$ is bounded locally in $C^{1,\ga}([T-1,T+1]\ti S^1)$ independently of $T$ and thus in 
$C^{2,\ga}([T-1,T+1]\ti S^1)$, for some $\ga\in (0,1)$. \smallskip

 \nind {\it 1-(ii) Convergence and sharp gradient estimates.}  
Let $\bar w(t)=\frac 1{2\gp}\int_0^{2\gp}w(t,\gth)d\gth$ be the average of $w(t,.)$ on $S^1$. Then 
$$-\bar w''(t)+ae^{(2-b\gg)t}H_0(t)+e^{(2-q)t}G_0(t)=0,
$$
where 
$$H_0(t)=\frac 1{2\gp}\int_{S^1}e^{bw}d\gth\quad\text{and }\; G_0(t)=-\frac m{2\gp}\int_{S^1}\left((\gg-w_t)^2+w_\gth^2\right)^\frac q2d\gth.
$$
Since $H_0$ and $G_0$ are bounded, $\bar w'(t)$ admits a limit when $t\to -\infty$, and this limit is necessarily equal to zero since $\bar w$ is bounded. Thus
$$\bar w'(t)+a\int_{-\infty}^te^{(2-b\gg)\gt}H_0(\gt)d\gt+\int_{-\infty}^te^{(2-q)\gt}G_0(\gt)d\gt=0.
$$
Therefore 
\bel{WT}|\bar w'(t)|\leq c\max\{e^{(2-b\gg)t}, e^{(2-q)t}\},
\ee
which implies that $\bar w(t)$ admits a finite limit $\ell$ when $t\to-\infty$. Set $w^*=w-\bar w$. Then 
\bel{W*}-w^*_{tt}-w^*_{\gth\gth}+ae^{(2-b\gg)t}H^*+e^{(2-q)t}G^*=0,
\ee
where 
$$H^*(t,.)=e^{bw}-\frac 1{2\gp}\int_{S^1}e^{bw}d\gth\quad\text{and }\; G^*(t)=\frac m{2\gp}\int_{S^1}\left|(\gg-w_t)^2+w_\gth^2\right|^\frac q2d\gth
-m\left(\gg-w_t)^2+w_\gth^2\right)^\frac q2.
$$
Since $\int_{S^1}w^*d\gth=0$, we have by Wirtinger inequality, 
$$-\int_{S^1} w_{\theta\theta}w^*d\gth=-\int_{S^1} w^*_{\theta\theta}w^*d\gth\geq \norm {w^*}^2_{L^2(S^1)},
$$
and
$$\int_{S^1}H^*(t,.)w^*d\gth=\int_{S^1}\left(e^{bw}-\frac 1{2\gp}\int_{S^1}e^{bw}d\gth\right)(w-\bar w)d\gs=
\int_{S^1}\left(e^{bw}-ae^{b\bar w}\right)(w-\bar w)d\gs\geq 0.
$$
Using the inequality
$$-\norm{w^*}_{L^2(S_1)}\frac{d^2}{dt^2} \norm {w^*}_{L^2(S^1)}\leq -\int_{S^1}w^*w^*_{tt}d\gth
$$
we obtain
$$-\frac{d^2}{dt^2} \norm {w^*}_{L^2(S^1)}+ \norm {w^*}_{L^2(S^1)}\leq L_0(t)e^{\gb t},
$$
where $\beta=\min\{2-q,2-b\gg\}>0$ and  $L_0(t)$ is some bounded function. Since $\norm {w^*}_{L^2(S^1)}$ is bounded and $\gb<1$, it implies
that 
 \begin{equation}\label{EX4}
\norm {w^*(t,.)}_{L^2(S^1)}\leq Ce^{\gb t}.
\end{equation}
Combining this estimate with the fact that $\bar w(t)\to \ell$ when $t\to-\infty$ we deduce that 
 \begin{equation}\label{EX5}
\lim_{t\to-\infty}\norm {w(t,.)-\ell}_{L^2(S^1)}=0.
\end{equation}
Let 
$$\CT[w]=\bigcup_{t\leq 0}\{(w(t,.),w_t)\}$$ 
be the negative trajectory of $(w,w_t)$. It is relatively compact in $C^2(S^1)\ti C^1(S^1)$ and because of the uniqueness of the limit,  we have that $w(t,.)\to \ell$ in $C^2(S^1)$ and thus $w_\gth(t,.)\to 0$ in $C^1(S^1)$. In order to 
make precise the behaviour of $w_t$, we differentiate $(\ref{W*})$ with respect to $t$ and get
\bel{W*T}-w^*_{ttt}-w^*_{t\gth\gth}+e^{(2-b\gg)t}\tilde H^*+e^{(2-q)t}\tilde G^*=0
\ee
where $\tilde H^*=a(2-b\gg)H^*-aH^*_t$ and $\tilde G^*=(2-q)G^*-G^*_t$ are bounded functions (from the $C^2$ a priori bounds on $w$). Hence by Wirtinger inequality again, 
$$-\frac{d^2}{dt^2} \norm {w_t^*}_{L^2(S^1)}+ \norm {w_t^*}_{L^2(S^1)}\leq \tilde L_0(t)e^{\gb t},
$$
where $\tilde L_0$ is a bounded function. Since $\norm {w_t^*}_{L^2(S^1)}$ is bounded on $(-\infty,0]$ it implies again that 
 \begin{equation}\label{EX6}
\norm {w_t^*(t,.)}_{L^2(S^1)}\leq C'e^{\gb t}.
\end{equation}
By $(\ref{WT})$, $\bar w_t(t)\to 0$ as $t\to-\infty$. This implies
 \begin{equation}\label{EX7}
\lim_{t\to-\infty}\norm {w_t(t,.)}_{L^2(S^1)}=0.
\end{equation}
As in the case of the function $w(t,.)$, the relative compactness of the set $\{w_t(t,.)\}_{t\leq 0}$ in $C^1(S^1)$ combined with $(\ref{EX7})$ implies that $w_t(t,.)\to 0$ in $C^1(S^1)$ as $t\to-\infty$. 
This implies $(\ref{EX2**})$.\smallskip


 \nind {\it 1-(iii) End of the proof.} 
Next we consider two solutions $u_\gg$ and $\tilde u_\gg$ satisfying $(\ref{Ex2})$ with $\gg\in ( 0,\frac 2b)$ and coinciding on $\prt B_1$. Denote $u=u_\gg$ and $\tilde u=\tilde u_\gg$. For $\epsilon\neq 0$ we recall that $E(x)=\ln\frac 1{|x|}$ and
 \begin{equation}\label {X3}U(r)=U_{\epsilon}(x)=u(x)+\epsilon E(x).
 \end{equation}
 Then
 $$\begin{array}{lll}\displaystyle-\Delta U+ae^{bU}-m|\nabla U|^q=ae^{bu}\left(e^{b\epsilon E}-1\right)-
 m\left[\left(|\nabla u|^2+\frac{\ge^2}{|x|^2}-2\frac{\ge}{|x|}\nabla u.{\bf i}\right)^\frac q2-|\nabla u|^q\right]\\[4mm]
 \phantom{+ae^{bU}-m|\nabla U|^q}
 \displaystyle=ae^{bu}\left(e^{b\epsilon E}-1\right)- m|\nabla u|^q\left[\left(1+\frac{\ge^2}{|x|^2|\nabla u|^2}-2\frac{\ge}{|x||\nabla u|^2}\nabla u.{\bf i}\right)^\frac q2-1\right].
 \end{array}$$
Then, when $r\to 0$, 
$$1+\frac{\ge^2}{|x|^2|\nabla u|^2}-2\frac{\ge}{|x||\nabla u|^2}\nabla u.{\bf i}=1+\frac{\ge^2}{\gg^2}(1+o(1))+2\frac{\ge }{\gg}(1+o(1)),
$$
thus
$$|\nabla u|^q\left[\left(1+\frac{\ge^2}{|x|^2|\nabla u|^2}-2\frac{\ge}{|x||\nabla u|^2}\nabla u.{\bf i}\right)^\frac q2-1\right]=\frac{\gg^q}{|x|^q}\left(\frac{q\ge}{\gg}+\frac{q(q-1)\ge^2}{2\gg}(1+o(1))\right).
$$
Since $b\gg\ln\frac1{|x|}+bK\geq bu(x)\geq b\gg\ln\frac1{|x|}-bK$, we have in the case where $\ge>0$:
$$ae^{bK}|x|^{-b\gg}(|x|^{-b\ge}-1)\geq ae^{bu}\left(e^{b\epsilon E}-1\right)\geq ae^{-bK}|x|^{-b\gg}(|x|^{-b\ge}-1),
$$
and in the case where $\ge<0$: 
$$ae^{bK}|x|^{-b\gg}(|x|^{-b\ge}-1)\leq ae^{bu}\left(e^{b\epsilon E}-1\right)\leq ae^{-bK}|x|^{-b\gg}(|x|^{-b\ge}-1).
$$ \smallskip

 \nind {\it 1-(iii)-(i): Case $\frac qb\leq \gamma<\frac 2b$}. we take $\ge>0$ and get
$$\CE(U)\geq ae^{-bK}|x|^{-b\gg}(|x|^{-b\ge}-1)-m\frac{\gg^q}{|x|^q}\left(\frac{q\ge}{\gg}+\frac{q(q-1)\ge^2}{2\gg}(1+o(1))\right).
$$
 Fix $r_0\in (0,1)$ and $\ge_0>0$ such that for any $0<r\leq r_0$ and $0<\ge\leq\ge_0$ there holds
$$m\gg^q\left(\frac{q\ge}{\gg}+\frac{q(q-1)\ge^2}{2\gg}(1+o(1))\right)\leq 2qm\gg^{q-1}\ge.
$$
Therefore 
$$\CE(U)\geq |x|^{-b\gg}\left(ae^{-bK}(r_0^{-b\ge}-1)-2qm\gg^{q-1}\ge|x|^{b\gg-q}\right).
$$
Up to reducing $r_0$ we have $\CE(U)\geq 0$ in $\overline B_{r_0}\setminus\{0\}$ for any $\ge\in (0,\ge_0]$.  Next, in $B_1\setminus B_{r_0}$ we have that 
$\nabla u_\gg$ is bounded, therefore
$$\CE(U)\geq -C\ge$$ 
for some constant $C>0$ whenever $\epsilon$ is small enough, Consequently this last inequality holds in whole $\overline B_1\setminus\{0\}$. Finally the function
$$x\mapsto U^*(x)=U(x)+\frac{1}{b}\ln (1+C\ge)
$$
is a supersolution of $(\ref{Ex1})$ larger than $\tilde u$ in a neighbourhood of $x=0$ and near $\prt B_1$. Hence $U^*\geq \tilde u_\gg$. Letting $\ge\to 0$ yields 
$u_\gg\geq \tilde u$. Similarly $\tilde u\geq  u$, which implies that  $\tilde u=  u$, that is $\tilde u_\gg=  u_\gg$.\smallskip

 \nind {\it 1-(iii)-(ii): Case $0\leq \gamma<\frac qb$}. We denote by $\eta$ the solution of $(\ref{eta})$
with value $\eta(x)=(2-q)^2(1-|x|^{2-q})$. For $\epsilon,\epsilon'>0$, set  $ U(x)=u(x)-\epsilon E(x)$ and $ U_1(x)= U(x)-\epsilon'\eta(x)$. Then
 $$\nabla U_1=\nabla U+(2-q)\ge'|x|^{1-q}{\bf i}\quad\text{and }\;-\Delta U_1=-\Delta U+\epsilon'\Delta \eta=-\Delta u-\epsilon'|x|^{-q},
 $$
 $$\BA{lll}\dsps|\nabla U_1|^q=\left(|\nabla U|^2+\frac{(2-q)^2\ge'^2}{|x|^{2(q-1)}}+\frac{2(2-q)\ge'}{|x|^{q-1}}\nabla U.{\bf i}\right)^\frac q2\\[4mm]
 \phantom{\dsps|\nabla U_1|^q}\dsps=|\nabla U|^q\left(1+\frac{(2-q)^2\ge'^2}{|x|^{2(q-1)}|\nabla U|^2}+\frac{2(2-q)\ge'}{|x|^{q-1}|\nabla U|^2}\nabla U.{\bf i}\right)^\frac q2,
 \EA$$
 $$e^{b U_1}=e^{b U}e^{-b\epsilon'\eta},
 $$
 and
 $$-\Gd U=-\Gd u=m|\nabla u|^q-ae^{b u}=m|\nabla ( U+\ge E)|^q-ae^{b( U+\ge E)}.
 $$
 Since $\eta\geq 0$ and $0<|x|\leq 1$, $e^{b U}\left(e^{-b\epsilon' \eta}-|x|^{-b\epsilon}\right)\leq 0$. Therefore
   \begin{equation}\label {X9} \begin{array}{lll}\displaystyle-\Delta  U_1+ae^{b U_1}-m|\nabla  U_1|^q=-\epsilon'|x|^{-q}+ae^{b U}\left(e^{-b\epsilon' \eta}-|x|^{-b\epsilon}\right) \\[2mm]\phantom{}\dsps
 -m|\nabla  U|^q\left[\left(1+\frac{(2-q)^2\ge'^2}{|x|^{2(q-1)}|\nabla  U|^2}+\frac{2(2-q)\ge'}{|x|^{q-1}|\nabla  U|^2}\nabla  U.{\bf i}\right)^\frac q2
-\left(1+\frac{\ge^2}{|x|^2|\nabla  U|^2}-\frac{2\ge}{|x||\nabla  U|^2}\nabla  U.{\bf i}\right)^\frac q2\right]\\
\\[0mm]\phantom{--}\dsps
\leq-\epsilon'|x|^{-q}
+C|\nabla  U|^{q}\left(\frac{\ge'}{|x|^{q-1}|\nabla  U|}+\frac{\ge'^2}{|x|^{2(q-1)}|\nabla  U|^2}+\frac{\ge}{|x||\nabla  U|}+\frac{\ge^2}{|x|^2|\nabla  U|^2}\right)(1+o(1)),
 \end{array}\end{equation}
for some $C>0$ as $r\to 0$.
In order $ U_1$ be a subsolution, the problem is reduced to proving that, under suitable conditions on $\epsilon$ and  $\epsilon'$, the expression
 $$A(x):=-\epsilon'|x|^{-q}
+2C|\nabla  U|^{q}\left(\frac{\ge'}{|x|^{q-1}|\nabla  U|}+\frac{\ge'^2}{|x|^{2(q-1)}|\nabla  U|^2}+\frac{\ge}{|x||\nabla  U|}+\frac{\ge^2}{|x|^2|\nabla  U|^2}\right)$$ 
is negative for $r$ small enough.
Because $\frac{2\gg} {|x|}\geq |\nabla  U(x)|\geq \frac\gg {2|x|}$ near $0$, there holds
 $$A(x)\leq -\epsilon'|x|^{-q}+2^{q+1}C\gg^q|x|^{-q}\left(\frac{2\ge'|x|^{2-q}}{\gg}+\frac{4\ge'^2|x|^{4-2q}}{\gg^2}+\frac{2\ge}{\gg}+\frac{4\ge^2}{\gg^2}\right).
 $$
We fix $\epsilon'=2^{q+1}C\gg^q\left(\frac{2\ge'}{\gg}+\frac{4\ge'^2}{\gg^2}+\frac{2\ge}{\gg}+\frac{4\ge^2}{\gg^2}\right)$, then 
 $$A(x)\leq 0\quad\text{for }\, 0<|x|\leq r_0.$$
Proceeding as in is the case  (3-i) there exists $\epsilon_0$ and $C>0$ such that, for any $0<\ge\leq \ge_0$ there holds 
\bel{EX1}\CE( U_1)\leq C\ge\quad\text{in }\; B_1\setminus\{0\}.
\ee 
The function $x\mapsto \tilde V_\ge(x):=\tilde u(x)+\frac{1}b\ln (1+C\ge)$ satisfies 
$$\CE(\tilde V_\ge)=C\ge e^{b\tilde u}\geq C.
$$
 Therefore $\tilde V_\ge$ is a super solution of $(\ref{EX1})$ which is larger than $ U_1$ near $x=0$ and $|x|=1$. Then by the maximum principle, $\tilde V_\ge\geq  U_1$. 
 Letting $\ge\to 0$ yields $\tilde u\geq u$. Similarly  $u\geq \tilde u$, thus $u_\gg=\tilde u_\gg$.\medskip

 
 \nind {\it 2- The case $\gg=\frac 2b$}.\smallskip
 
 \nind {\it 2-(i) Gradient estimates}. We first show that 
 \bel{EX2a}
 \lim_{x\to 0}\left(u(x)-\frac 2b\left(\ln\frac{1}{|x|}-\ln(1-\ln |x|)\right)\right)=\ell\quad\text{and }\,\lim_{x\to 0}|x|\nabla u(x)=-\frac 2b{\bf i}
 \ee
 for some $\ell\in [-K,K]$ with $K$ from inequality $(\ref{AE-8})$ in \rcor{glob} and ${\bf i}=\frac{x}{|x|}$.\smallskip
 
 We set again $w(t,.)=u(r,.)+\frac 2b\ln r$ with $t=\ln r$. Thus $w$ satisfies
$$-\frac {2}{b}\ln (1-t)-K\leq w(t,.)\leq -\frac {2}{b}\ln (1-t)+K
$$
and
\bel{EQ++}w_{tt}+w_{\theta\theta}-ae^{bw}+me^{(2-q)t}\left((\tfrac 2b-w_t)^2+w_\gth^2\right)^\frac q2=0.
\ee
Introducing again $z(t,.)=(\ln (1-t))^{-1}w(t,.)$, then 
$$ -\frac 2b-\frac K{\ln (1-t)}\leq z(t,.)\leq -\frac 2b+\frac K{\ln (1-t)}
$$
and $z$ which is bounded satisfies equation $(\ref{sp3})$; thus $w_t$ satisfies $(\ref{sp3+})$. Since $z$ is bounded and negative and all the coefficients are smooth and bounded, the standard regularity theory implies that $z$ is bounded in $C^{2,\ga}([T'-1,T'+1]\ti S^1)$ independently of $T'$, thus 
\bel{sp6*}
\lim_{t\to-\infty}z(t,.)=-\frac 2b\quad\text{in }\, C^2(S^1)
\ee
Set
$$\gl(t,.)=w(t,.)+\frac{2}{b}\ln(1-t)
$$
Then $\gl$ is uniformly bounded and it satisfies in $(-\infty,0]\ti S^1$,
\bel{sp6-1}-\gl_{tt}-\gl_{\gth\gth}-\frac{2}{b(1-t)^2}+\frac{a}{(1-t)^2}e^{b\gl}-me^{(2-q)t}\left((\tfrac 2b-w_t)^2+w_\gth^2\right)^\frac q2=0.
\ee
Let $\bar\gl(t)$ be the circular average  of $\gl(t,.)$ and $\gl^*(t,.)=\gl(t,.)-\bar\gl(t)$. Then 
$$-\gl^*_{tt}-\gl^*_{\gth\gth}+\frac{a}{(1-t)^2}\left(e^{b\gl}-\overline{e^{b\gl}}\right) +e^{(2-q-\ge)t}G^*_\ge=0
$$ 
where $G^*_\ge$ is a bounded function. Again, by Wirtinger inequality we obtain that 
$$
\norm{\gl^*(t,.)}_{L^2(S^1)}\leq C\max\left\{e^{t},e^{(2-q-\ge)t}\right\}=Ce^{(2-q-\ge)t},
$$
hence
\bel{sp6-3}\lim_{t\to-\infty}\norm{\gl^*(t,.)}_{L^2(S^1)}=0
\ee
By Gagliardo-Nirenberg inequality we have for any $p>1$
\bel{sp6-4}\BA{lll}\norm{\nabla\gl^*(t,.)}_{W^{1,p}(S^1)}\leq C\norm{\nabla\gl^*(t,.)}^{1-\frac 1p}_{W^{2,2}(S^1)}\norm{\nabla\gl^*(t,.)}^{\frac 1p}_{L^{2}(S^1)}\\
[2mm]\phantom{\norm{\nabla\gl^*(t,.)}_{W^{1,p}(S^1)}}\leq 
C\norm{\nabla\gl^*(t,.)}^{1-\frac 1p}_{C^{2}(S^1)}\norm{\nabla\gl^*(t,.)}^{\frac 1p}_{L^{2}(S^1)}.
\EA\ee
Since $\gl$ is bounded in $C^2(S^1)$, we combined $(\ref{sp6-4})$ with Morrey-Sobolev inequality and obtain, by taking $p>2$, 
\bel{sp6-5}
\norm{\nabla\gl^*(t,.)}_{C^{0,1-\frac 1p}(S^1)}\leq \norm{\gl^*(t,.)}^{\frac 1p}_{L^{2}(S^1)}\leq Ce^{\gd t},
\ee
where $\gd=\frac{2-q-\ge}{p}>0$. The function $\bar\gl$ satisfies 
\bel{sp6-12}-\bar\gl_{tt}-\frac{2}{b(1-t)^2}+\frac{a}{(1-t)^2}\overline{e^{b\gl}} +e^{(2-q-\ge)t}G_\ge=0
\ee
If it is monotone at $-\infty$, then it converges to some $\ell\in [-K,K]$. If it oscillates infinitely many times when $t\to-\infty$, there exists a sequence $\{t_n\}$ tending to 
$-\infty$ such that $\bar\gl (t)$ achieves a local maximum when $t=t_{2n}$ and a local minimum when $t=t_{2n+1}$. Hence 
$$\bar\gl_{tt}(t_{2n})\leq 0 \,\text{ and }\;\bar\gl_{tt}(t_{2n+1})\geq 0.
$$
We have $$e^{b\gl(t,.)}=e^{b\bar\gl(t)}(e^{\gl(t,.)-\bar\gl(t)})=e^{b\bar\gl(t)}\left(1+D(t,.)e^{\delta t}\right)$$
for some bounded function $D$, hence 
$$\overline{e^{b\gl}}=e^{b\bar\gl(t)}\left(1+\overline{D(t,.)}e^{\delta t}\right)
$$
Plugging this estimate into $(\ref{sp6-12})$ we obtain that there exists a bounded function $B(t)$ and $\gn>0$ such that 
\bel{sp6-13}-\frac{2}{b}+a e^{b\bar\gl(t_n)} +e^{\gn t_n}B(t_n)=(1-t_n)^2\bar\gl_{tt}(t_n)
\ee
This yields
\bel{sp6-14}
e^{b\bar\gl(t_{2n})}\leq \frac{2}{ab}-e^{\gn t_{2n}}B(t_{2n})\quad\text{and }\;e^{b\bar\gl(t_{2n+1})}\geq \frac{2}{ab}-e^{\gn t_{2n+1}}B(t_{2n+1})
\ee
and finally
\bel{sp6-15}
\limsup_{t\to-\infty}\bar\gl(t)=\liminf_{t\to-\infty}\bar\gl(t)=\lim_{t\to-\infty}\bar\gl(t)=\frac 1b\ln\frac{2}{ab}:=\ell.
\ee
This implies 
\bel{sp6-16}
 \lim_{x\to 0}\left(u(x)-\frac 2b\left(\ln\frac{1}{|x|}-\ln(1-\ln |x|)\right)\right)=\ell.
 \ee
 
To end the estimate we have still to prove that $w_t\to 0$. Differentiating $(\ref{EQ++})$ with respect to $t$ we obtain 
\bel{sp7}
-w_{ttt}-w_{t\gth\gth}+abe^{bw}w_t+e^{(2-q)t}\tilde G=0
\ee
where $\tilde G=O(\ln (1-t))$. Since $e^{bw}\geq \frac{e^{-bK}}{(1-t)^2}$ we get
$$-\frac{d^2}{dt^2}\norm{w_t}_{L^2(S^1)}+\frac{abe^{-bK}}{(1-t)^2}\norm{w_t}_{L^2(S^1)}\leq Ae^{\gm t}
$$
where $\gm<2-q$ and $A>0$. Put $X(t)=\norm{w_t(t,.)}_{L^2(S^1)}+\frac{A'}{\gm^2}e^{\gm t}$ with $A'>A$. Then there exists $T^*<0$ such that for $t\leq T^*$  we have
$$
-X''+\frac {e^{-bK} ab}{(1-t)^{2}}X\leq 0.
$$
The corresponding equation
$$
\CL(x):=-x''+\frac {e^{-bK} ab}{(1-t)^{2}}x=0\quad\text{with } t\leq 0
$$ 
admits two linearly independent solutions
$$x_1(t)=(1-t)^{\gl_1}\quad\text{ and }\; x_2(t)=(1-t)^{\gl_2},
$$
when $t<0$ arbitrary and $\gl_1<0<\gl_2$. Since $X(t)=o(x_2(t))$ when $t\to-\infty$, it follows that $X(t)=O(x_1(t))$. Hence 
\bel{sp6}
\norm {w_t}_{L^2(S^1)}\leq B(1-t)^{\gl_1}\quad\text{for }t\leq 0
\ee
for some $B\geq 0$. In order to have a uniform estimate, we use Poisson representation formula for the solution $U$ of the Dirichlet problem for the Laplacian in the disk of radius R
\bel{sp9*}
U(r,\gth)=\frac{R^2-r^2}{2\gp R}\int_{0}^{2\gp}\frac{U(R,\gf))d\gf}{\sqrt{r^2+R^2-2rR\cos(\gth-\gf)}}
\ee
Using the variable $t=\ln r$ and $\tilde U(t,\gth)=U(r,\gth)$, it follows that the bounded solution of 
\bel{sp10*}\BA{lll}
-\tilde U_{tt}-\tilde U_{\gth\gth}=0\qquad &\text{in }(-\infty,T)\ti S^1\\
\phantom{,,,,}
\tilde U(T,\gth)=\tilde U_T(\gth)&\text{in } S^1,
\EA\ee
is expressed by
\bel{sp11}
\tilde U(t,\gth)=\frac{e^{2T}-e^{2t}}{2\gp e^T}\int_{0}^{2\gp}\frac{\tilde U(T,\gf))d\gf}{\sqrt{e^{2t}+e^{2T}-2e^{t+T}\cos(\gth-\gf)}}
\ee
Hence, for $T=\frac t2$ 
\bel{sp12}
\tilde U(t,\gth)\leq \frac{1+e^{\frac t2}}{2\gp}\int_{0}^{2\gp}\tilde U_+(\tfrac t2,\gf)d\gf
\ee
From $(\ref{sp7})$  there exists some $ A^*\geq 0$ such that the function $W(t,\gth)=w_t(t,\gth)+A^*e^{\gm t}$ satisfies 
$$-W_+\,_{tt}-W_+\,_{\gth\gth}\leq 0\qquad\text{in }(-\infty,0)\ti S^1.
$$
Since 
$$\norm{W(T,.)}_{L^1(S^1)}=\norm{w_t(T,.)+A^*e^{\gm T}}_{L^1(S^1)}\leq 2\gp A^*e^{\mu T}+\sqrt{2\gp}\norm{w_t(T,.)}_{L^2(S^1)},
$$
we derive from $(\ref{sp6})$ and $(\ref{sp12})$, that there holds for some $C>0$,
\bel{sp13}
(w_t(t,\gth)+A^*e^{\gm t})_+=W_+(t,\gth)\leq C((t+1)^{\gl_1}\qquad\text{for all }(t,\gth)\in (-\infty,0]\ti S^1.
\ee
This implies
\bel{sp14}
\lim_{x\to 0}|x|\nabla u(x)=-\frac 2b{\bf i}.
\ee
The end of  the proof is similar to the case $\frac qb<\gg<\frac 2b$.
\qeda
\medskip



 \mysection{Singular solutions in the case $q>2$}
 In all this section we assume $q>2$ and $m,b,a>0$.
 \subsection{The eikonal equation}
 The eikonal equation associated to $(\ref{I-1})$ 
 is 
 \begin{equation}\label{E1}
ae^{bw}=m|\nabla w|^q\quad\text{in } \BBR^2\setminus\{0\}.
\end{equation}
This equation is completely integrable since it is equivalent to
$$\frac bq\left(\frac am\right)^{\frac 1q}=|\nabla(e^{-\frac bqw})|.
$$
By \cite{CaCr} $w$ is radial with respect to $0$.
Since any solution cannot have extremum, it is either decreasing or increasing with respect $|x|$. If it is decreasing on $(0,1]$ with $u(0)=c$, then 
\begin{equation}\label{E3}
w(x):=w_c(x)=\frac qb\ln\left(\frac{qm^{\frac 1q}}{ba^\frac 1q|x|+qm^{\frac 1q}e^{-\frac{bc}{q}}}\right)
\end{equation}
Letting $c\to\infty$ yields
\begin{equation}\label{E4}
\lim_{c\to\infty}w_c(x):=w_\infty(x)=\frac qb\ln\left(\frac{qm^{\frac 1q}}{ba^{\frac 1q}|x|}\right)=\frac qb\ln\left(\frac{1}{|x|}\right)+\frac qb\ln\left(\frac{qm^{\frac 1q}}{ba^{\frac 1q}}\right).
\end{equation}
This is also a solution of
\begin{equation}\label{E5}
-\Delta u+ae^{bu}-m|\nabla u|^q=0
\end{equation}
in the whole punctured plane $\BBR^2\setminus\{0\}$. It vanishes on $\partial B_1$ if and only if $qm^{\frac 1q}=ba^{\frac 1q}$. \medskip
\subsection{Singular solutions} 
 \blemma{grad2} If $u$ is a radial nonnegative solution of $(\ref{E5})$ in $B_1\setminus\{0\}$ it is decreasing on some $(0,r^*]\subset (0,1]$ and there holds
  \begin{equation}\label{E23-5}
-u'(r)=|u'(r)|\geq cr^{-\frac {1}{q-1}}\qquad\text{for all  }r\in (0,r^*],
\end{equation}
for some $c>0$,
 \es
\Proof  If $u$ is a positive radial function satisfying $(\ref{E5})$ in $B_1\setminus\{0\}$, it cannot have any maximum there in, hence either it is decreasing on $(0,1)$ which implies that $u(0)<u(1)$, or it is decreasing in some interval $(0,r^*)$ and increasing in $(r^*,1)$. For simplicity and without loss of generality, we will assume that 
 $r^*=1$ (this is always the case if $u(1)=0$).  We have
 $$-u''-\frac 1ru'-m|u'|^q\leq 0,
 $$
 then with $V=-u'$, we obtain 
 $$V'+\frac 1rV-mV^q\leq 0\Longleftrightarrow V'V^{-q}+\frac 1rV^{1-q}-m\leq 0\Longleftrightarrow (r^{1-q}V^{1-q})' +m(q-1)r^{1-q}\geq 0$$
 Therefore the function
 $$r\mapsto r^{1-q}V^{1-q}(r) -m\frac{q-1}{q-2}r^{2-q}
 $$
 is increasing, which implies
  \begin{equation}\label{E14*}
V^{1-q}(r) \leq m\frac{q-1}{q-2}r(1+o(1))\quad\text{as }r\to 0.
\end{equation}
Replacing $V$ by $-u'$ we obtain 
  \begin{equation}\label{E15}
u'(r)\leq -\left(\frac{q-2}{m(q-1)}\right)^{\frac{1}{q-1}}r^{-\frac{1}{q-1}}(1+o(1))\quad\text{as }r\to 0.
\end{equation}
\qeda
\bprop{Proposition 7}  Let $u$ be a positive solution of $(\ref{E5})$ in $B_1\setminus\{0\}$. The following statements are equivalent:\smallskip

\noindent (i) $e^{bu}\in L^1(B_1)$\smallskip

\noindent (ii) $|\nabla u|^q\in L^1(B_1)$. \smallskip

\noindent (iii) $u$ can be extended to in $\overline B_1$ as a H\"older continuous function and 
the extended function, still denoted by $u$, is a solution of
   \begin{equation}\label{E23-6}-\Delta u+ae^{bu}-m|\nabla u|^q=0 \qquad \text{in } {\bf D}'( B_1).
\end{equation}
\es
\Proof {\it $(i)\Longrightarrow (ii)$}. If  $e ^{bu}\in L^1(B_1)$, it follows by Brezis-Lions Lemma \cite{BrLi} that $|\nabla u|^q\in L^1(B_1)$ and there exists $\gamma\geq 0$ such that 
   \begin{equation}\label{E23-7}
   -\Delta u+ae^{bu}-m|\nabla u|^q=2\pi\gamma\delta_0 \qquad \text{in } {\bf D}'( B_1).
\end{equation}
 The average function $\bar u$ satisfies
\begin{equation}\label{E12}-(r\bar u_r)_r=r\left(m\overline{(u_r^2+r^{-2}u_\gth^2)^{\frac q2}}-a\overline{e^{bu}}\right).
\end{equation}
 Therefore $(r\bar u_r)'\in L^1(0,1)$ which implies that $r\bar u_r(r)$ admits a finite limit when $r\to 0$. If this limit is not zero, say $\tilde\gamma$, then 
 $$\bar u_r(r)=\frac{\tilde\gamma}{r}(1+o(1))\Longrightarrow |\bar u_r(r)|^q\sim \frac{|\tilde\gamma|^q}{r^q} . 
 $$ 
 Now 
 $$\BA{lll}\dsps |u_r(r)|\leq \frac 1{2\gp}\int_{S^1}\left|u_r(r,.)\right|d\gth\leq \frac{1}{2\gp}\int_{S^1}|\nabla u(r,.)|d\gth
 \leq \left(\frac{1}{2\gp}\int_{S^1}\left|\nabla u(r,.)\right|^qd\gth\right)^{\frac1q}
 \EA
 $$
 Since $q>2$ we have $|\nabla u|^q\notin L^1(B_1)$, contradiction. Hence $\tilde \gg=0$.
Therefore $\bar u_r(r)=o(r^{-1})$, and this implies (\ref{E13}).

 \noindent {\it $(ii)\Longrightarrow (i)$}. We proceed by contradiction and assume that  $e ^{bu}\notin L^1(B_1)$.  If $|\nabla u|^q\in L^1(B_1)$, then 
 $$
r\bar u_r(r)-\bar u_r(1)=\int_r^1\left(m\overline{|\nabla u|^q}-a\overline{e^{bu}}\right)tdt=-a\int_r^1\overline{e^{bu}}tdt+O(1)\to-\infty\quad\text{as }r\to 0.
 $$
 Hence $ r\bar u_r(r)\to-\infty$ when $r\to 0$, which would imply
 $$\lim_{r\to 0}\frac{u(r)}{-\ln r}=\infty,
 $$
 which would contradict the a priori estimate. Therefore, if  $ e ^{bu}$ is not integrable in $B_1$  implies that it is the same holds with $|\nabla u|^q$. This proves the claim. \smallskip
 
\noindent {\it $(iii)\Longrightarrow (i)$}. Since $u$ is continuous in $\overline B_1$, $e ^{bu}\in L^1(B_1)$.\smallskip
 
\noindent {\it $(ii)\Longrightarrow (iii)$}.  Since $|\nabla u|\in L^{q}(B_1)$ and $u$ is $C^1$ in $\overline B_1\setminus\{0\}$, it follows by Morrey-Sobolev theorem that 
it can be extended by continuity to $x=0$ and that the resulting function, still denoted by $u$ belongs to $C^{0,1-\frac 2q}(\overline B_1)$. 
In particular 
\begin{equation}\label{E13}|u(x)-u(0)|\leq c|x|^{1-\frac 2q}.
\end{equation}
Furthermore, by Brezis-Lions lemma, there exists $\gg\geq 0$ such that $(\ref{E23-7})$ holds. By the same averaging method that in the case $(i)\Longrightarrow (ii)$, we have also
$$\lim_{r\to 0}\frac{\bar u(r)}{-\ln r}=\gg.
$$
The previous estimate implies that $\gg=0$. This ends the proof. 
\qeda \medskip

\noindent\Remark If $u$ is a radial function estimate $(\ref{E13})$ can be improved. Actually, since by \rlemma{grad2} one has $(\ref{E23-5})$, the integration of this inequality yields
\begin{equation}\label{E14}|u(r)-u(0)|\leq cr^{\frac {q-2}{q-1}},
\end{equation}
and $\frac {q-2}{q-1}>\frac {q-2}{q}$.
\subsection{Proof of Theorem 3}

The proof of Theorem 3 is as follows: we first observe that the conclusion of Theorem 3 holds in the class of radial solutions. Then we prove that if $u$ is a positive solution $(\ref{E5})$ in $B_1\setminus\{0\}$ which tends to 
$\infty$ when $x\to 0$, then there exists a positive solution $\tilde u$ which satisfies 
$$\left(u-\max\{u(y):|y|=1\}\right)_+\leq\tilde u\leq u\quad\text{in }B_1\setminus\{0\}
$$
and vanishes on $\prt B_1$. This part of the proof is delicate and based upon a combination of Bernstein method for estimating the gradient of a solution and a subtle variation of 
Boccardo-Murat-Puel classical existence result of solution in presence of a couple of ordered supersolution and subsolution. Finally we prove by adaptating the moving plane method that the function $\tilde u$ is radial and we apply the  radial result and the above estimate.\smallskip

The following result has already been obtained in the proof of assertion 2-(ii) in \rth{class}.

\bprop {Proposition 8}  Let $q>2$ and $u$ be a positive solution of $(\ref{E5})$ in $B_1\setminus\{0\}$. Then 
   \begin{equation}\label{E23-9}
 \int_{B_1}|\nabla u|^q dx=\infty \Longrightarrow  \liminf_{x\to 0}\frac{u(x)}{-\ln |x|}\geq\frac 2b.
\end{equation}
\es
\smallskip

\medskip

\noindent\Remark Using the estimates of Proposition 2 we have precisely
   \begin{equation}\label{E23-14}
u(x)\geq \frac 2b \left(\ln\frac 1{|x|}-\ln(1-\ln |x|)\right) -K.
\end{equation}\medskip

Under the assumption of \rprop{Proposition 8} a radial solution $u$ vanishing on $\partial B_1$ is decreasing. In the next result we prove that $ru'(r)$ is monotone. \medskip

\blemma { Lemma 2} Let $u$ be a positive radial singular solution of $(\ref{E5})$ in $B_1\setminus\{0\}$, then  $ru'(r)$ is monotone near $r=0$.\es
\smallskip

 \noindent {\it Proof}. We set
    \begin{equation}\label{E23-17*}
\CW(r)=ae^{bu(r)}-m|u'(r)|^q.
\end{equation}
Using the equation $(\ref{E5})$ we have
$$\CW'=u'ae^{bu}\left(b-qm|u'|^{q-2}\right)+  qm|u'|^{q-2}u'\left(|u'|^{q}+\frac{u'}{r}\right).
$$
Replacing $ae^{bu}$ by $\CW+m|u'|^q$, we obtain
$$\CW'=u'\left(b-qm|u'|^{q-2}\right)\CW+mb|u'|^qu'-(m-1)qm|u'|^{2q-2}u'+\frac{qm}{r}|u'|^q
$$
Denote
$$\Theta(r)=u'(r)\left(b-qm|u'|^{q-2}(r)\right),
$$
and 
$$\Sigma(r)=mb|u'|^qu'-(m-1)qm|u'|^{2q-2}u'+\frac{qm}{r}|u'|^q
$$
By $(\ref{E23-5})$, $u'(r)$ is nonpositive near $r=0$ and tends to $-\infty$ when $r\to 0$. Then
$$\Sigma(r)=(m-1)q|u'|^{2q-1}(1+o(1))+\frac{qm}{r}|u'|^q\quad\text{as }r\to 0,
$$
since $q>2$. By integration on $[r,r_0]\subset (0,1)$, 
$$\frac{d}{dr}\left(e^{\int_r^{r_0}\Theta(s)ds}\CW(r)\right)=e^{\int_r^{r_0}\Theta(s)ds}\Sigma (r),
$$
hence
$$\CW(r)=e^{-\int_r^{r_0}\Theta(s)ds}\CW(r_0)-\int_r^{r_0}e^{-\int_r^{s}\Theta(t)dt}\Sigma (s)ds.
$$
There exists $r^*\in (0,1)$ such that $\Sigma(r)>0$ for $(0<r\leq r^*$.  Therefore if for some $r_0\in (0,r^*]$ there holds $\CW(r_0)\leq 0$, it follows that 
$\CW(r)< 0$ for any $r\in (0,r_0)$. Since
$$(ru')'=r\left(ae^{bu}-m|u'|^q\right)=r\CW,
$$
this implies that $(ru')'\leq 0$ on $(0,r_0]$. On the contrary, for any if for any $r\in (0,r^*]$
we have $\CW(r)>0$, we conclude that $(ru')'\leq 0$ on $(0,r_0]$.\qeda\medskip

\bprop {Proposition 9} Let $q>2$ and $u$ be a positive radial solution of $(\ref{E5})$ in $B_1\setminus\{0\}$ satisfying  $\int_{B_1}|\nabla u|^q dx=\infty$. If 
$ru'(r)$ is monotone then 
   \begin{equation}\label{E23-15}
 \lim_{x\to 0}\frac{u(x)}{-\ln |x|}=\frac qb.
\end{equation}
\es
\smallskip

 \noindent \Proof If $u$ satisfies the above assumption, there holds by Propositions 2 and 8,  
    \begin{equation}\label{E23-16}
\frac 2b\ln\frac 1{|x|}-K(b)\ln(1-\ln |x|)\leq u(x)\leq \frac qb\ln\left(\frac{1}{|x|}\right)+\frac qb\left[\left(\ln\left(\frac b{qm^{\frac 1q}}\right)\right)_+-\ln\left(\frac b{qm^{\frac 1q}}\right) \right],
\end{equation}
where $K(b)=1$ if $b\geq 2$ and $K(b)=\frac 2b$ if $0<b\leq 2$. We set $t=\ln r$ and $w(t)=u(r)$. Then 
   \begin{equation}\label{E23-17}-w_{tt}+ae^{bw+2t}-me^{(2-q)t}|w_t|^q=0.
\end{equation}
Moreover, if $ru'(r)$ is monotone, then $w_t$ is monotone too. From $(\ref{E23-16})$,  
$$
-\frac 2b t-K(b)\ln(1-t)\leq w(t)\leq -\frac qb t+\frac qb\left[\left(\ln\left(\frac b{qm^{\frac 1q}}\right)\right)_+ -\ln\left(\frac b{qm^{\frac 1q}}\right)\right]\quad\text{on }(-\infty,0).
$$
If $w_t$ is monotone, it admits a limit $\ell$, finite or not, when $t\to -\infty$. The only possibility is that 
$$
-\frac qb\leq \ell\leq -\frac 2b.
$$
and $w(t)=\ell t(1+o(1))$ when $t\to-\infty$. From $(\ref{E23-17})$,
$$w_{tt}=-me^{(2-q)t}|\ell|^q(1+o(1))+ae^{(b\ell+2)t+o(t)}.
$$
There always holds $2-q\geq b\ell+2$. In order that $w_{tt}$ does not tend to $-\infty$ (which would be a contradiction), necessarily $2-q=b\ell+2$ which implies $\ell=-\frac qb$. This yields $(\ref{E23-15})$.\qeda
\medskip

In next result we obtain an a priori estimate of the gradient of a solution by a combination of the Bernstein and the Keller-Osserman methods as in  
\cite[Lemma 3.1]{BV} (see also \cite[Lemma 2.2]{BVGHV*} for a simpler version).
\bprop {Bern} If $q>2$ and $u$ is a positive solution of of $(\ref{E5})$ in $B_1\setminus\{0\}$, then for any $\gd>0$ there exists $C^*=C^*(a,b,m,q,\gd)>0$ such that 
      \begin{equation}\label{E25-1}
|\nabla u(x)|\leq \frac{C^*}{|x|}\quad\text{for all }x\in B_{1-\gd}\setminus\{0\}.
\ee
\es
\Proof If $u$ is such a solution, we set 
$$z=|\nabla u|^2.
$$
Then by an explicit computation and using Bochner identity and Schwarz inequality
      \begin{equation}\label{E25-2*}-\frac 12\Gd z=-|D^2u|^2 -\nabla(\Gd u).\nabla u\leq -\frac 12(\Gd u)^2-\nabla(\Gd u).\nabla u.
\ee
Since $\Gd u=ae^{bu}-m z^{\frac q2}$, we have for any $\ge>0$ and explicit $C(\ge)>0$,
$$\BA{lll}\dsps -\Gd z\leq -(ae^{bu}-m z^{\frac q2})^2-2\nabla(ae^{bu}-m z^{\frac q2}).\nabla u\\[1mm]
\phantom{-\Gd z}\dsps\leq -a^2e^{2bu}-m^2z^q+2ame^{bu}z^{\frac q2}-2abe^{bu}z+qmz^{\frac q2}\frac{|\nabla z|}{\sqrt z}
\\[1mm]
\phantom{-\Gd z}\dsps\leq -m^2z^q +2am\left(\ge z^q+C(\ge)e^{2bu}\right)+qm\left(\ge z^q+C(\ge)\frac{|\nabla z|^2}{z}\right).
\EA$$
Therefore
$$\BA{lll}\dsps
-\Gd z+\left(m^2-\ge\left(2am+qm\right)\right)z^q\leq 2amC(\ge)e^{2bu}+qmC(\ge)\frac{|\nabla z|^2}{z}.
\EA$$
We choose $\ge$ such that $m^2-\ge\left(2am+qm\right)\geq \frac {m^2}2$. By the a priori estimate $(\ref{I-12})$,  $e^{2bu}\leq e^{2b\tilde \Gl}|x|^{-2q}$
Therefore $z$ satisfies 
   \begin{equation}\label{E25-2}
   -\Gd z+\frac {m^2}2z^q\leq \frac{C_1}{|x|^{2q}}+C_2\frac{|\nabla z|^2}{z},
   \end{equation}
   for some $C_j=C_j(a,b,m,q)>0$, $j=1,2$. By \cite[Lemma 3.1]{BV},  it follows that 
      \begin{equation}\label{E25-3}
z(x)\leq C^*_1|x|^{-\frac{2}{q-1}}+C^*_2|x|^{-2}\leq (C_1+C_2)|x|^{-2}.
   \end{equation}
This implies $(\ref{E25-1})$.\qeda\medskip

\noindent\Remark When $1<q<2$, the above method leads to the estimate
      \begin{equation}\label{E25-4}
|\nabla u(x)|\leq \frac{C^*}{|x|^{-\frac 1{q-1}}}\quad\text{for all }x\in B_{1-\gd}\setminus\{0\},
\ee
which is not optimal for radial functions as \rlemma{grad1} shows it.

\blemma {exist} Let $q>2$ and $u\in C^1(B_1\setminus\{0\})\cap C(\overline{B_1}\setminus\{0\})$ be a nonnegative solutions of $(\ref{E5})$ in $B_1\setminus\{0\}$ which satisfies 
   \begin{equation}\label{E26-1*}\liminf_{x\to 0}\frac{u(x)}{-\ln |x|}\geq\frac 2b.
   \end{equation}
   If $\gm=\max\{u(x):|x|=1\}$, there exists a positive solution $\tilde u\in C^1(B_1\setminus\{0\})\cap C(\overline{B_1}\setminus\{0\})$ of $(\ref{E5})$ in $B_1\setminus\{0\}$ which vanishes on $\prt B_1$ and satisfies
      \begin{equation}\label{E26-2}(u-\gm)_+\leq \tilde u\leq u\quad\text{in }B_1\setminus\{0\}.
   \end{equation}
\es
\Proof {\it Step 1: Construction of an approximate solution}. The function $u_\gm=(u-\gm)_+$ is a positive subsolution of $(\ref{E5})$ in $B_1\setminus\{0\}$ which vanishes on $\prt B_1$ and it is dominated by the solution $u$. However we cannot apply directly the classical result of Boccardo, Murat and Puel \cite{BMP} since $q>2$ (the limit is $q=2$) and a variation of the construction is needed. 
For $\ge\in (0,1)$ it follows by \rprop {Bern} that 
   \begin{equation}\label{E26-2+}\BA{lll}\dsps
|\nabla u(x)|\leq \frac{C^*}{\ge}\qquad\text{for all }x\in B_1\setminus B_\ge.
\EA
\ee
and $C^*$ is independent of $\ge$. Up to changing $C^*$ the above inequality holds if $u$ is replaced by $(u-\gm)_+$.
We define now the function $\phi_\ge$ on $\BBR_+$ by 
   \begin{equation}\label{E26-3}
\phi_\ge(\xi)=\left\{\BA{lll}\xi^{\frac q2}\quad&\text{if }0\leq \xi\leq  \left( \frac{C^*}{\ge}\right)^2\\[2mm]\dsps
\frac{q(q-1)}{2}\left(\frac {C^*}\ge\right)^{q-2}\xi+q(2-q)\left(\frac {C^*}\ge\right)^{q-1}\sqrt \xi+\frac{q^2-3q+2}{2}\left(\frac {C^*}\ge\right)^{q}\quad&\text{if } \xi\geq  \left( \frac{C^*}{\ge}\right)^2.
\EA\right.
\ee
The function $\xi\mapsto \phi_\ge\left(|\xi|^2\right)$ belongs to $C^2(\BBR^2)$ and is convex with a quadratic growth at infinity. The result from \cite{BMP} applies. Since $\phi_\ge\left(|\nabla u|^2\right)=|\nabla u|^q$ (resp. 
$\phi_\ge\left(|\nabla (u-\gm)_+|^2\right)=|\nabla (u-\gm)_+|^q$), the problem 
   \begin{equation}\label{E26-4}\BA{lll}
-\Gd w+ae^{bw}-m\phi_\ge\left(|\nabla w|^2\right)=0\qquad&\text{in }B_1\setminus\overline B_\ge\\[1mm]
\phantom{-\Gd +ae^{bw}-m\phi_\ge\left(|\nabla w|^2\right)}
w=0\qquad&\text{on }\prt B_1\\[1mm]
\phantom{-\Gd +ae^{bw}-m\phi_\ge\left(|\nabla w|^2\right)}
w=(u-\gm)_+&\text{on }\prt B_\ge,
\EA
\end{equation}
admits a solution $w_\ge$ which satisfies 
   \begin{equation}\label{E26-4*}
(u-\gm)_+\leq w_\ge\leq u\qquad\text{in }B_1\setminus\overline B_\ge.
\ee
{\it Step 2: Local estimate of $\nabla w_\ge$}. 
In order to have a local estimate independent of $\ge$, we apply Bernstein method as in \rprop{Bern}. Setting again $z=|\nabla w_\ge|^2$ we have from the equation 
$(\ref{E26-4})$,
$$\BA{lll}-\Gd z\leq \left(ae^{bw_\ge}-m\phi_\ge\left(z\right)\right)^2-2\nabla \left(ae^{bw_\ge}-m\phi_\ge\left(z\right)\right).\nabla w_\ge\\[2mm]
\phantom{-\Gd z}\dsps
\leq -a^2e^{2bw_\ge}-m^2\phi^2_\ge\left(z\right)+2ame^{bw_\ge}\phi_\ge\left(z\right)-2abe^{bw_\ge}z+2m\phi'_\ge\left(z\right)|\nabla z|\sqrt z
\\[2mm]
\phantom{-\Gd z}\dsps
\leq -(m^2-2am\gd^2)\phi^2_\ge\left(z\right)+\frac{2am}{\gd^2}e^{2bw_\ge}+2m\gd^2\phi'^2_\ge(z)z^2+\frac{2m}{\gd^2}\frac{|\nabla z|^2}{z}.
\EA$$
By the definition of $\phi_\ge$ there exists $\gk>0$ independent of $\ge$ such that $\phi'^2_\ge(z)z^2\leq \gk\phi^2_\ge\left(z\right)$. Provided $\gd>0$ is small enough and using $(\ref{AE-6})$ and using estimate $(\ref{E26-4})$ we obtain,
\bel{E26-5}
-\Gd z+\frac{m^2}2\phi^2_\ge\left(z\right)\leq \frac{2am}{\gd^2}e^{2bw_\ge}+\frac{2m}{\gd^2}\frac{|\nabla z|^2}{z}
\leq \frac{K}{|x|^{2q}}+\frac{2m}{\gd^2}\frac{|\nabla z|^2}{z}.
\ee
By the expression of $\phi_\ge(\xi)$, there exists constants $\gth$ and $C_\gth>0$ such that 
\bel{E26-6}
\phi^2_\ge(z)\geq 2\gth z^2-C_\gth
\ee
Therefore $(\ref{E26-5})$ combined with $(\ref{E26-4})$ and $(\ref{I-12})$ yields 
\bel{E26-7}
-\Gd z+ m^2\gth z^2\leq \frac{K}{|x|^{2q}}+\frac{C_\gth m^2}2+\frac{2m}{\gd^2}\frac{|\nabla z|^2}{z}\quad\text{in }B_1\setminus B_\ge.
\ee
If $\ge<|x|<1-\ge$ we apply \cite[Lemma 3.1]{BV}  in $B_{|x|-\ge}(x)$ with $q=2$ and derive 
\bel{E26-8}
z(x):=|\nabla w_\ge(x)|^2\leq \frac {\tilde C_1}{ |x|^{q}}+\frac{\tilde C_2}{(|x|-\ge)^2}
\ee
 with $\tilde C_j>0$ independent of $\ge$. We now define $\ge^*>\ge$ by the identity
 \bel{E26-9}
\frac {\tilde C_1}{ \ge^{*\,q}}+\frac{\tilde C_2}{(\ge^*-\ge)^2}=\left(\frac{C^*}{\ge}\right)^2
\ee
where $C^*$ is defined at $(\ref{E26-3})$ (such a $\ge^*$ exists by the mean value theorem). The value of $\ge^*$ cannot be explicited but 
 \bel{E26-10}
\tilde C_1^{\frac 1q}\left(\frac {\ge}{C^*}\right)^{\frac 2{q}}\leq \ge^*\leq \tilde C_1^{\frac 1q}\left(\frac {\ge}{C^*}\right)^{\frac 2{q}}(1+o(1))\quad\text{as }\ge\to 0.
\ee
As a consequence the function $w_\ge$ is also a solution of  
   \begin{equation}\label{E26-11}\BA{lll}
-\Gd w+ae^{bw}-m|\nabla w|^q=0\qquad&\text{in }B_1\setminus\overline B_{\ge^*}\\[1mm]
\phantom{-\Gd +ae^{bw}-m|\nabla w|^q}
w=0\qquad&\text{on }\prt B_1,
\EA
\ee
and thus it satisfies the upper estimate obtained in \rprop{Bern},
      \begin{equation}\label{E26-12}
|\nabla w_\ge(x)|\leq \frac{C^*}{|x|-\ge^*}\quad\text{for all }x\in B_{1}\setminus\overline B_{\ge^*}.
\ee
\medskip

\nind {\it Step 3: End of the proof.} It follows from $(\ref{E26-12})$ that for any $\ge\in (0,1)$ small enough, $w_\ge$ and $\nabla w_\ge $ are uniformly bounded in 
$B_{1}\setminus\overline B_{2\ge^*}$. By the standard elliptic equations regularity theory, $w_\ge$ remains bounded in  $C^{2,\ga}(K)$ for any $K\subset B_{1}\setminus\overline B_{3\ge^*}$. Using a diagonal sequence, there exists $\{\ge_n\}$ converging to $0$ and $\tilde u\in C^2(\overline B_1\setminus \{0\})$ such that $\{w_{\ge_n}\}$ converges to 
$\tilde u$ in the $C^2_{loc}$-topology of $\overline B_1\setminus \{0\}$. Therefore $\tilde u$ is a solution of $B_1\setminus\{0\}$ which vanishes on $\prt B_1$ and satisfies $(u-\gm)_+\leq \tilde u\leq u$ jointly with the gradient estimate $(\ref{E26-2+})$.\qeda


\bprop {symm} Let $q>2$ and $u\in C^1(B_1\setminus\{0\})\cap C(\overline{B_1}\setminus\{0\})$ is a nonnegative solutions of $(\ref{E5})$ in $B_1\setminus\{0\}$ which vanishes on 
$\prt B_1$ and satisfies 
   \begin{equation}\label{E26-1}
 \lim_{x\to 0}u(x)=\infty.
\end{equation}
Then $u$ is radially symmetric.
\es
\Proof The proof is an easy variant of the moving plane method as it is introduced in \cite{GNN} (see also \cite[Theorem 8.2.2]{PS}). The only difference is the assumption that $u$ satisfies $(\ref{E26-1})$ which actually allows the reflection through the singular point. The  sketch of the proof is as follows. We set $x=(x_1,x_2)$ and, if $\gl\in (0,1)$, $B^\gl=\{x\in B_1:x_1>\gl\}$, $T^\gl=\{x\in B_1:x_1=\gl\}$ and 
$C^\gl=\prt B_1\cap \overline{B^\gl}$. Put $x_\gl=(2\gl-x_1,x_2)$
and $u_\gl(x)=u(x_\gl)$. If $w=u_\gl-u$,  there holds
   \begin{equation}\label{E26-2++}-\Gd w(x)+A(x)w(x)-B(x).\nabla w(x)=0\quad\text{in }B^\gl
\end{equation}
where 
$$A(x)=ab\int_0^1e^{b(tu_\gl(x)+(1-t)u(x)}dt
$$
and 
$$B(x)=mq\int_0^1|\nabla \left(tu_\gl(x)+(1-t)u(x)\right)|^{q-2}\nabla \left(tu_\gl(x)+(1-t)u(x)\right)dt
$$
The two functions $A$ and $B$ are H\"older continuous in $B^\gl$ for $\gl>0$.  Since $u>0$ in $B_1\setminus\{0\}$ and vanishes on  $\prt B_1$, the outward normal derivative  $\frac{\prt u}{\prt{\bf n}}$ is negative on $\prt B_1$. For $\gl\in (\frac12,1)$ the function $w$ satisfies $(\ref{E26-2++})$ in $B^\gl$. It vanishes on $T^\gl$ and is positive on $C^\gl$ except at the two end-points 
belonging to $\prt B_1\cap T^\gl$. In fact, if we assume now that there exists $x^*\in \overline {B^\gl}$ where $w$ achieves a negative minimum, then $x^*\in B^\gl$ and because $w$ is $C^2$, 
$\nabla w(x^*)=0$ and $-\Gd w(x^*)\leq 0$. Since 
$$-\Gd w(x^*)+A(x^*)w(x^*)-B(x^*).\nabla w(x^*)=0
$$
and $A(x^*)>0$ we have a contradiction. If $\gl\in (0,\frac 12]$ the function $w$ is not continuous in whole $\overline {B_\gl}$, but since $u$ satisfies $(\ref{E26-1})$, the function 
$w$ tends to infinity when $x\to (2\gl,0)$ (this point is the reflected image of $(0,0)$ by the symmetry $(x\mapsto x_\gl)$). Hence the minimum of $w$ in $ \overline {B_\gl}$ coincides with the minimum in $ \overline {B^\gl}\setminus B_\gd((2\gl,0))$. Since $w$ is $C^2$ therein, this implies that it cannot achieve a negative minimum. As a consequence $w\geq 0$ in 
$B^0_+:=B_1\cap\{x=(x_1,x_2):x_1>0\}$, and by continuity this inequality holds also in $\overline B_1\cap\{x=(x_1,x_2):x_1\geq 0\}\setminus\{(0,0)\}$. This implies in particular 
$u_{x_1}(0,x_2)\geq 0$. Starting the reflection process from the point $(-1,0)$ we obtain in the same way $u_{x_1}(0,x_2)\leq 0$. Hence $u_{x_1}(0,x_2)=0 $ for all $x_2\in (-1,1)$. Because the equation is invariant by rotation we can replace 
the axis $(0,x_2)$ by any axis through $(0,0)$. Consequently we have $u_\gth(\gth,r)=0$ in polar coordinates, which is the claim.\qeda\medskip

\nind {\it Proof of Theorem 3}. It follows from \rprop{Proposition 7}, \rprop{Proposition 8}, \rprop{Proposition 9}, \rlemma{exist} and \rprop{symm}.\qeda

\medskip

\nind\Remark This type of dichotomy which occurs in this case between regular solutions and only one type of singular solutions is well known for semilinear Lane-Emden 
equations
   \begin{equation}\label{E23-20}
   -\Delta u=u^p\quad\text{in } B_1\setminus\{0\}
\end{equation}
when $N\geq 3$ and $p>\frac{N}{N-2}$. To our knowledge it is a novelty for equation of $(\ref{E5})$  and in a 2-dimensional framework. 

\end{document}